\documentclass[12pt]{amsart}
\usepackage{pgfplots}
\usepackage{tikz}

\usetikzlibrary{automata}
\usetikzlibrary{positioning}
\usepackage{amsmath, amssymb} 
\usepackage{amsthm}
\usepackage{graphicx}
\usepackage[shortlabels]{enumitem}

\usepackage{amsfonts,geometry,hyperref}

\newtheorem{conj}{Conjecture}[section]

\newtheorem{theorem}{Theorem}[section]

\newtheorem{cor}[theorem]{Corollary}
\newtheorem{question}[theorem]{Question}

\allowdisplaybreaks

\begin{document}

\title[Hansen-Vuki{\v{c}}evi{'c} Conjecture]{Commuting Conjugacy Class Graphs of Finite Groups and the Hansen-Vuki{\v{c}}evi{'c} Conjecture}

\author[S. Das, R. K. Nath, Y. Shang]{Shrabani Das,  Rajat Kanti Nath* and Yilun Shang*}

\address{S. Das, Department of Mathematical Sciences, Tezpur University, Napaam-784028, Sonitpur, Assam, India.\newline
Department of Mathematics, Sibsagar University, Joysagar-785665, Sibsagar, Assam, India.
}

\email{shrabanidas904@gmail.com}

\address{R. K. Nath, Department of Mathematical Sciences, Tezpur University, Napaam-784028, Sonitpur, Assam, India.} 
\email{ rajatkantinath@yahoo.com}

\address{Y. Shang, Department of Computer and Information Sciences, Northumbria University, Newcastle NE1 8ST, UK.
} 
\email{yilun.shang@northumbria.ac.uk}
\thanks{*Correspondence: rajatkantinath@yahoo.com (R. K. Nath); yilun.shang@northumbria.ac.uk (Y. Shang)}

\begin{abstract}
In this work, we compute the first and second Zagreb indices for the commuting conjugacy class graphs associated with finite groups. We identify multiple classes of finite groups whose commuting conjugacy class graphs are shown to satisfy the Hansen-Vuki{\v{c}}evi{\'c} conjecture. Specifically, we prove that the conjecture holds for the commuting conjugacy class graphs of dihedral groups ($D_{2m}$), dicyclic groups, semidihedral groups, and various other two-generator groups. Moreover, we examine the case where the quotient $G/Z(G)$ is isomorphic to $D_{2m}$, $\mathbb{Z}_p \times \mathbb{Z}_p$, a Frobenius group of order $pq$ or $p^2q$, or any group of order $p^3$, for primes $p$ and $q$. In each of these cases, we demonstrate that the corresponding commuting conjugacy class graph satisfies the Hansen-Vuki{\v{c}}evi{\'c} conjecture.
\end{abstract}

\thanks{ }

\maketitle
\section{Introduction}

Let $\mathcal{G}$ denote the collection of all graphs. A topological index for a graph $\Gamma \in \mathcal{G}$ is expressed as $T(\Gamma)$, where $T$ is a real-valued function mapping from $\mathcal{G}$ to $\mathbb{R}$. Originally, such indices were introduced to characterize various chemical attributes of molecular structures. Since 1947, numerous topological indices have been formulated based on different graph parameters. The earliest of these is the Wiener index, a distance-based measure introduced by Wiener \cite{Wiener-97}. Within the category of degree-based indices, the Zagreb indices are historically the first and were proposed by Gutman and Trinajsti{\'c} in 1972 \cite{Gut-Trin-72}. These indices were initially employed to investigate how the total $\pi$-electron energy of a molecule correlates with its structural features. As mentioned in \cite{Z-index-30y-2003}, the Zagreb indices have also been applied in the analysis of molecular complexity, chirality, ZE-isomerism, and heterosystems, among others. Over time, mathematicians have explored the general mathematical characteristics of a variety of topological indices. A comprehensive overview of the mathematical aspects of Zagreb indices is available in the survey article \cite{Gut-Das-2004}. More recent developments on Zagreb indices can be found in e.g. \cite{ismail, noureen, raza}.

Let $\Gamma$ be a simple, undirected graph characterized by the vertex set $v(\Gamma)$ and edge set $e(\Gamma)$. Two well-known degree-based topological indices for such graphs are the first and second Zagreb indices, denoted respectively by $M_{1}(\Gamma)$ and $M_{2}(\Gamma)$. These are formally defined as follows:
\[
M_{1}(\Gamma) = \sum\limits_{v \in v(\Gamma)} \deg(v)^{2} \quad \text{and} \quad M_{2}(\Gamma) = \sum\limits_{uv \in e(\Gamma)} \deg(u)\deg(v),
\]
where $\deg(v)$ denotes the degree of a vertex $v$, that is, the number of edges incident to it. In an effort to explore the relationship between these two indices, Hansen and Vuki{\v{c}}evi{\'c} \cite{hansen2007comparing} proposed a conjecture in 2007, which compares the behavior of $M_1(\Gamma)$ and $M_2(\Gamma)$ for various classes of graphs.

\begin{conj}\label{Conj}
(Hansen-Vuki{\v{c}}evi{\'c} Conjecture) Assume that $\Gamma$ is a simple finite graph. We have the following inequality. 
\begin{equation}\label{Conj-eq}
 \dfrac{M_{2}(\Gamma)}{\vert e(\Gamma) \vert} \geq \dfrac{M_{1}(\Gamma)}{\vert v(\Gamma) \vert} .
\end{equation}
\end{conj}

It was demonstrated in \cite{hansen2007comparing} that the conjecture does not hold in general, as evidenced by the counterexample $\Gamma = K_{1, 5} \sqcup K_3$. Nevertheless, Hansen and Vuki{\v{c}}evi{\'c} \cite{hansen2007comparing} verified that Conjecture~\ref{Conj} remains valid within the class of chemical graphs. The validity of the conjecture has also been established for trees in \cite{vukicevic2007comparing}, with equality in equation~\eqref{Conj-eq} occurring specifically when $\Gamma$ is a star graph. For connected unicyclic graphs, the conjecture was confirmed in \cite{liu2008conjecture}, where equality is attained in the case of cycle graphs. A more detailed examination of equality conditions in~\eqref{Conj-eq} was conducted in \cite{vukicevic2011some}, and a broader survey on the comparison of Zagreb indices is available in \cite{Liu-You-2011}.

Despite these advances, the complete classification of graphs that either satisfy or violate Conjecture~\ref{Conj} remains an open area of investigation. Recently, Das et al.~\cite{SD-AS-RKN-2023} identified various finite non-abelian groups whose commuting or non-commuting graphs fulfill the Hansen-Vuki{\v{c}}evi{\'c} conjecture. Similarly, Das and Nath \cite{SD-RKN-2024} provided several examples of finite groups whose super commuting graphs also conform to the conjecture.

In this paper, we focus on the commuting conjugacy class graphs of finite groups. Specifically, we compute their first and second Zagreb indices and identify several families of finite groups for which these graphs satisfy the Hansen-Vuki{\v{c}}evi{\'c} conjecture. Section~2.1 is devoted to analyzing the commuting conjugacy class graphs of certain group families. These include the dihedral groups $D_{2m}$ defined by the presentation $\langle a, b : a^m = b^2 = 1, bab^{-1} = a^{-1} \rangle$ for $m \geq 3$; the dicyclic groups $Q_{4m}$ presented as $\langle a, b : a^{2m} = 1, b^2 = a^m, bab^{-1} = a^{-1} \rangle$ for $m \geq 2$; and the semidihedral groups $SD_{8m}$ with presentation $\langle a, b : a^{4m} = b^2 = 1, bab = a^{2m-1} \rangle$ for $m \geq 2$.

In addition to these, we examine the group $V_{8m}$ given by $\langle a, b : a^{2m}=b^{4}=1, ba=b^{-1}a^{-1}, b^{-1}a=a^{-1}b \rangle$ for $m \geq 1$; the group $U_{(n,m)}$ defined by $\langle a, b : a^{2n}=b^m=1, a^{-1}ba=b^{-1} \rangle$ for $m \geq 3$, $n \geq 2$; and the group $G(p,m,n)$ with the presentation $\langle a,b : a^{p^m}=b^{p^m}=[a,b]^p=1, [a,[a,b]]=[b,[a,b]]=1 \rangle$ for any prime $p$ and integers $m, n \geq 1$. Section~2.2 explores commuting conjugacy class graphs of groups $G$ for which the quotient group $G/Z(G)$ is isomorphic to one of several specific types. These include $D_{2m}$, $\mathbb{Z}_p \times \mathbb{Z}_p$, Frobenius groups of order $pq$ and $p^2q$, and all groups of order $p^3$, where $p$ and $q$ are prime numbers.

The commuting conjugacy class graph of a group $G$, abbreviated as the CCC-graph, is defined as a simple undirected graph in which each vertex corresponds to a conjugacy class of a non-central element of $G$. Two distinct vertices, denoted $x^G$ and $y^G$—where $x^G$ represents the conjugacy class containing an element $x \in G$—are joined by an edge if and only if there exist representatives $x' \in x^G$ and $y' \in y^G$ such that $x'y' = y'x'$. This graph is denoted by $\mathcal{CCC}(G)$ and was first introduced by Herzog et al.~\cite{MH-ML-MM-2009} in 2009. Subsequent research has explored various structural properties of $\mathcal{CCC}(G)$. In 2016, Mohammadian et al.~\cite{AM-AE-DGMF-BW-2016} classified all finite groups $G$ for which the CCC-graph $\mathcal{CCC}(G)$ is triangle-free. Further developments were made in \cite{MAS-ARA-2020, MAS-ARA-2020-2}, where Salahshour and Ashrafi studied the structures of $\mathcal{CCC}(G)$ for several families of finite CA-groups. Additional work by Salahshour in \cite{MAS-2020} examined CCC-graphs in cases where the quotient group $G/Z(G)$ is isomorphic to a dihedral group. More recently, Rezaei et al.~\cite{MR-ZF-2024} determined the structure of $\mathcal{CCC}(G)$ when $G/Z(G)$ corresponds to a Frobenius group of order $pq$ or $p^2q$, where $p$ and $q$ are distinct prime numbers. Spectral characteristics of CCC-graphs have also been explored in a series of papers, including \cite{BN-2021, BN-2024, JN-2025, JNP-2024}. In \cite{BN-2022}, the genus of $\mathcal{CCC}(G)$ was computed, and various families of groups were identified whose CCC-graphs are planar, toroidal, double-toroidal, or even triple-toroidal. A comprehensive overview of the CCC-graph and its applications to group theory is available in the survey article \cite{CJNS-2024}.

\section{Zagreb indices of CCC-graphs}

The following result, obtained in \cite{SD-AS-RKN-2023}, is useful in our computations.

\begin{theorem}\label{thm1}
Suppose that $\Gamma$ is the disjoint union of the graphs $\Gamma_{1}, \Gamma_{2}, \dots, \Gamma_{n}$. Let $\Gamma_{i} = l_{i}K_{m_{i}}$ for $ i = 1, 2, \dots, k$, where $K_{m_{i}}$'s are complete graphs on $ m_{i} $ vertices and $l_iK_{m_i}$ is the disjoint union of $l_i$ copies of $K_{m_i}$. We obtain 
\[
M_{1}(\Gamma) = \sum_{i=1}^{k}l_{i}m_{i}(m_{i}-1)^{2} ~~~ \text{and} ~~~ M_{2}(\Gamma) = \sum_{i=1}^{k}l_{i}\dfrac{m_{i}(m_{i}-1)^{3}}{2}.
\]
\end{theorem}

\subsection{Zagreb indices of CCC-graphs of certain well-known finite groups}

In this section, we consider the CCC-graphs for the groups $D_{2m}, Q_{4m}, V_{8m}, SD_{8m}, U_{(n,m)}$ and $G(p,m,n)$ and compute their Zagreb indices to finally check if they satisfy Hansen-Vuki{\v{c}}evi{\'c} conjecture. The structures of the CCC-graphs of the above mentioned groups are obtained in \cite{MAS-ARA-2020-2}. However, while obtaining the structure of the CCC-graph of $U_{(n,m)}$ in \cite{MAS-ARA-2020-2}, one case was missed, which has been pointed out by Jannat and Nath in \cite{JN-2024}.

\begin{theorem}
If $G$ is the dihedral group of order $2m$, where $m \geq 3$, then 
    \[
    M_1(\mathcal{CCC}(G)) = \begin{cases}
        \frac{(m-1)(m-3)^2}{8}, & ~~~~\text{ when } m \text{ is odd} \\
        \frac{(m-2)(m-4)^2}{8}, & ~~~~\text{ when } m \text{ is even and } \frac{m}{2} \text{ is even} \\
        \frac{(m-2)(m-4)^2}{8} +2, & ~~~~\text{ when } m \text{ is even and } \frac{m}{2} \text{ is odd}
    \end{cases}
    \]
    and
    \[
    M_2(\mathcal{CCC}(G)) = \begin{cases}
        \frac{(m-1)(m-3)^3}{32}, & ~~~~\text{ when } m \text{ is odd} \\
        \frac{(m-2)(m-4)^3}{32}, & ~~~~\text{ when } m \text{ is even and } \frac{m}{2} \text{ is even} \\
        \frac{(m-2)(m-4)^3}{32} +1, & ~~~~\text{ when } m \text{ is even and } \frac{m}{2} \text{ is odd.}
    \end{cases}
    \]
    Further, $\frac{M_{2}(\mathcal{CCC}(G))}{|e(\mathcal{CCC}(G))|} \geq \frac{M_{1}(\mathcal{CCC}(G))}{|v(\mathcal{CCC}(G))|}$ with equality when $m=3, 4$ and $6$.
\end{theorem}
\begin{proof}
     \textbf{Case 1.} $ m $ is odd.

   We have $\mathcal{CCC}(G)=K_{\frac{m-1}{2}} \sqcup K_1$. Therefore, by Theorem \ref{thm1}, we get the expressions for $M_1(\mathcal{CCC}(G))$ and $M_2(\mathcal{CCC}(G))$.
   
   Note that $|v(\mathcal{CCC}(G))|=\frac{m-1}{2}+1=\frac{m+1}{2}$ and $|e(\mathcal{CCC}(G))|=\binom{\frac{m-1}{2}}{2}=\frac{(m-1)(m-3)}{8}$. Therefore,
\[
\frac{M_{1}(\mathcal{CCC}(G))}{|v(\mathcal{CCC}(G))|}=\frac{(m-1)(m-3)^2}{4(m+1)} \quad \text{ and } \quad \frac{M_{2}(\mathcal{CCC}(G))}{|e(\mathcal{CCC}(G))|}=\frac{(m-3)^2}{4}.
\]
Hence,
\[
\frac{M_{2}(\mathcal{CCC}(G))}{|e(\mathcal{CCC}(G))|} - \frac{M_{1}(\mathcal{CCC}(G))}{|v(\mathcal{CCC}(G))|} 
= \frac{(m-3)^2}{4}\left( \frac{2}{m+1}\right) 
 \geq 0 
\]
for all $m \geq 3$.

\vspace{.2cm}

\noindent     \textbf{Case 2.} $ m $ is even and $\frac{m}{2}$ is even.

  We have $\mathcal{CCC}(G)=K_{\frac{m}{2}-1} \sqcup 2K_1$. Therefore, by Theorem \ref{thm1}, we 
  get the expressions for $M_1(\mathcal{CCC}(G))$ and $M_2(\mathcal{CCC}(G))$. 
  
  Note that $|v(\mathcal{CCC}(G))|=\frac{m}{2}-1+2=\frac{m+2}{2}$ and $|e(\mathcal{CCC}(G))|=\binom{\frac{m}{2}-1}{2}=\frac{(m-2)(m-4)}{8}$. 
Therefore,
\[
\frac{M_{1}(\mathcal{CCC}(G))}{|v(\mathcal{CCC}(G))|}=\frac{(m-2)(m-4)^2}{4(m+2)} \quad \text{ and } \quad \frac{M_{2}(\mathcal{CCC}(G))}{|e(\mathcal{CCC}(G))|}=\frac{(m-4)^2}{4}.
\]
Hence,
\[
\frac{M_{2}(\mathcal{CCC}(G))}{|e(\mathcal{CCC}(G))|} - \frac{M_{1}(\mathcal{CCC}(G))}{|v(\mathcal{CCC}(G))|} 
= \frac{(m-4)^2}{4}\left( \frac{4}{m+2}\right) 
 \geq 0 
\]
for all $m \geq 4$.

\vspace{.2cm}

\noindent \textbf{Case 3.} $m$ is even and $\frac{m}{2}$ is odd.

 We have $\mathcal{CCC}(G)=K_{\frac{m}{2}-1} \sqcup K_2$. Therefore, by Theorem \ref{thm1}, we get the expressions for $M_1(\mathcal{CCC}(G))$ and $M_2(\mathcal{CCC}(G))$. 
 
 Note that $|v(\mathcal{CCC}(G))|=\frac{m}{2}-1+2=\frac{m+2}{2}$ and $|e(\mathcal{CCC}(G))|=\binom{\frac{m}{2}-1}{2}+1=\frac{(m-2)(m-4)}{8}+1$. 
Therefore,
\[
\frac{M_{2}(\mathcal{CCC}(G))}{|e(\mathcal{CCC}(G))|} = \frac{(m-2)(m-4)^3 + 32}{4((m-2)(m-4) + 8)}  \text{ and } \frac{M_{1}(\mathcal{CCC}(G))}{|v(\mathcal{CCC}(G))|} = \frac{(m-2)(m-4)^2 + 16}{4(m+2)}.
\]
Hence,
\begin{align*}
\frac{M_{2}(\mathcal{CCC}(G))}{|e(\mathcal{CCC}(G))|} &- \frac{M_{1}(\mathcal{CCC}(G))}{|v(\mathcal{CCC}(G))|}\\ 
&=\frac{m^2(m^2-16m+68)+20m(m-10)+8(m+18)}{(m+2)(m-2)(m-4)+8(m+2)}
:=\frac{f(m)}{g(m)}.
\end{align*}
Note that  $g(m)=16\times|e(\mathcal{CCC}(G))| \times |v(\mathcal{CCC}(G))|>0$. We have $f(6)=0$. Also, since $m(m-16)+68>0$ for all $m \geq 10$, therefore $\frac{f(m)}{g(m)} \geq 0$. Hence, the result follows.
\end{proof}
 \begin{theorem}
If $G$ is the dicyclic group of order $4m$, where $m \geq 2$, then      
    \begin{align*}
     M_1(\mathcal{CCC}(G)) = \begin{cases}
        (m-1)(m-2)^2, & \text{ when } m \text{ is even} \\
        (m-1)(m-2)^2+2, & \text{ when } m \text{ is odd}
    \end{cases}
    \end{align*}
    and
    \begin{align*}
    M_2(\mathcal{CCC}(G)) = \begin{cases}
        \frac{(m-1)(m-2)^3}{2}, & \text{ when } m \text{ is even} \\
        \frac{(m-1)(m-2)^3}{2}+1, & \text{ when } m \text{ is odd.}
        \end{cases}
    \end{align*}
    Further, $\frac{M_{2}(\mathcal{CCC}(G))}{|e(\mathcal{CCC}(G))|} \geq \frac{M_{1}(\mathcal{CCC}(G))}{|v(\mathcal{CCC}(G))|}$ with equality when $m=2$ and $3$.
\end{theorem}
\begin{proof}
    \textbf{Case 1.} $m$ is even.

     We have $\mathcal{CCC}(G)=K_{m-1} \sqcup 2K_1$. Therefore, by Theorem \ref{thm1}, we 
     get the expressions for $M_1(\mathcal{CCC}(G))$ and $M_2(\mathcal{CCC}(G))$.
     
     Note that $|v(\mathcal{CCC}(G))|=m-1+2=m+1$ and $|e(\mathcal{CCC}(G))|=\binom{m-1}{2}=\frac{(m-1)(m-2)}{2}$. Therefore,
\[
\frac{M_{2}(\mathcal{CCC}(G))}{|e(\mathcal{CCC}(G))|} = (m-2)^2 \text{ and } \frac{M_{1}(\mathcal{CCC}(G))}{|v(\mathcal{CCC}(G))|} = \frac{(m-1)(m-2)^2}{m+1}.
\]
Hence,
\[
\frac{M_{2}(\mathcal{CCC}(G))}{|e(\mathcal{CCC}(G))|} - \frac{M_{1}(\mathcal{CCC}(G))}{|v(\mathcal{CCC}(G))|} 
=\frac{2(m-2)^2}{m+1} 
 \geq 0 
\]
for all  $m \geq 2$.

\vspace{.2cm}

\noindent \textbf{Case 2.} $m$ is odd.

 We have $\mathcal{CCC}(G)=K_{m-1} \sqcup K_2$. Therefore, by Theorem \ref{thm1}, we get the expressions for $M_1(\mathcal{CCC}(G))$ and $M_2(\mathcal{CCC}(G))$.
 
 Note that $|v(\mathcal{CCC}(G))|=m-1+2=m+1$ and $|e(\mathcal{CCC}(G))|=\binom{m-1}{2}+1=\frac{(m-1)(m-2)}{2}+1$. Therefore,
 \[
\frac{M_{2}(\mathcal{CCC}(G))}{|e(\mathcal{CCC}(G))|} = \frac{(m-1)(m-2)^3+2}{(m-1)(m-2)+2} \text{ and } \frac{M_{1}(\mathcal{CCC}(G))}{|v(\mathcal{CCC}(G))|} = \frac{(m-1)(m-2)^2+2}{m+1}.
\]
Hence,
\[
\frac{M_{2}(\mathcal{CCC}(G))}{|e(\mathcal{CCC}(G))|} - \frac{M_{1}(\mathcal{CCC}(G))}{|v(\mathcal{CCC}(G))|} 
=\frac{m^2(m^2-8m+14)+8m(m-3)+9}{(m+1)(m-1)(m-2)+2(m+1)} 
:=\frac{f(m)}{g(m)}.
\]
Note that  $g(m)=2 \times|e(\mathcal{CCC}(G))| \times |v(\mathcal{CCC}(G))|>0$. We have $f(3)=0$ and $f(5)=64 >0$. Since $m(m-8)+14 > 0$ for all $m \geq 7$, we have $\frac{f(m)}{g(m)} \geq 0$. Hence the result follows.
\end{proof}
\begin{theorem}\label{CCC SD_8m}
	If $G$ is the semidihedral group of order $8m$, where $m \geq 2$, then 
	
	\begin{align*}
		M_1(\mathcal{CCC}(G)) = \begin{cases}
			(2m-1)(2m-2)^2, & \text{ when } m \text{ is even} \\
			(2m-2)(2m-3)^2+36, & \text{ when } m \text{ is odd}
		\end{cases}
	\end{align*}
	and
	\begin{align*}
		M_2(\mathcal{CCC}(G)) = \begin{cases}
			\frac{(2m-1)(2m-2)^3}{2}, & \text{ when } m \text{ is even} \\
			(m-1)(2m-3)^3+54, & \text{ when } m \text{ is odd.}
		\end{cases}
	\end{align*}
	Further, $\frac{M_{2}(\mathcal{CCC}(G))}{|e(\mathcal{CCC}(G))|} \geq \frac{M_{1}(\mathcal{CCC}(G))}{|v(\mathcal{CCC}(G))|}$ with equality when $m=3$.
\end{theorem}
\begin{proof}
	\textbf{Case 1.} $m$ is even.
	
	We have $\mathcal{CCC}(G)=K_{2m-1} \sqcup 2K_1$. Therefore, by Theorem \ref{thm1}, we get the expressions for $M_1(\mathcal{CCC}(G))$ and $M_2(\mathcal{CCC}(G))$.
	
	 Note that $|v(\mathcal{CCC}(G))|=2m-1+2=2m+1$ and $|e(\mathcal{CCC}(G))|=\binom{2m-1}{2}+0=\frac{(2m-1)(2m-2)}{2}$. Therefore,
\[
\frac{M_{2}(\mathcal{CCC}(G))}{|e(\mathcal{CCC}(G))|} = (2m-2)^2 \text{ and } \frac{M_{1}(\mathcal{CCC}(G))}{|v(\mathcal{CCC}(G))|} = \frac{(2m-1)(2m-2)^2}{2m+1}.
\]
Hence,
\[
\frac{M_{2}(\mathcal{CCC}(G))}{|e(\mathcal{CCC}(G))|} - \frac{M_{1}(\mathcal{CCC}(G))}{|v(\mathcal{CCC}(G))|} 
=\frac{2(2m-2)^2}{2m+1} 
\geq 0 
\]
for all  $m \geq 1$. 

\vspace{.2cm}
	
	\noindent \textbf{Case 2.} $m$ is odd.
	
	We have $\mathcal{CCC}(G)=K_{2m-2}\sqcup K_4$. Therefore, by Theorem \ref{thm1}, we get the expressions for $M_1(\mathcal{CCC}(G))$ and $M_2(\mathcal{CCC}(G))$.
	
	Note that  $|v(\mathcal{CCC}(G))|=2m-2+4=2m+2$ and $|e(\mathcal{CCC}(G))|=\binom{2m-2}{2}+\binom{4}{2}=\frac{(2m-2)(2m-3)}{2}+6=(m-1)(2m-3)+6$. 
Therefore,
\[
\frac{M_{2}(\mathcal{CCC}(G))}{|e(\mathcal{CCC}(G))|} = \frac{(m-1)(2m-3)^3+54}{(m-1)(2m-3)+6} \text{ and } \frac{M_{1}(\mathcal{CCC}(G))}{|v(\mathcal{CCC}(G))|} = \frac{(m-1)(2m-3)^2+18}{m+1}.
\]
Hence,
\begin{align*}
		\frac{M_{2}(\mathcal{CCC}(G))}{|e(\mathcal{CCC}(G))|} - \frac{M_{1}(\mathcal{CCC}(G))}{|v(\mathcal{CCC}(G))|} 
		&=\frac{8m^3(m-7)+4m(30m-18)}{(m+1)(m-1)(2m-3)+6(m+1)} 
		:=\frac{f(m)}{g(m)}.
	\end{align*}
	Note that $g(m)=\frac{1}{2} \times |e(\mathcal{CCC}(G))| \times |v(\mathcal{CCC}(G))|>0$. We have $f(3)=0$ and $f(5)=640>0$. Also, since $30m-18>0$ for all $m \geq 7$, therefore $\frac{f(m)}{g(m)} \geq 0$. Hence the result follows.
\end{proof}

\begin{theorem}\label{CCC V_8m}
If $G = V_{8m}$, where $m \geq 1$, then 
     \begin{align*}
     M_1(\mathcal{CCC}(G)) = \begin{cases}
        (2m-2)(2m-3)^2+4, & \text{ when } m \text{ is even} \\
        (2m-1)(2m-2)^2, & \text{ when } m \text{ is odd}
    \end{cases}
    \end{align*}
    and
    \begin{align*}
    M_2(\mathcal{CCC}(G)) = \begin{cases}
        (m-1)(2m-3)^3+2, & \text{ when } m \text{ is even} \\
        \frac{(2m-1)(2m-2)^3}{2}+1, & \text{ when } m \text{ is odd.}
        \end{cases}
    \end{align*}
    Further, $\frac{M_{2}(\mathcal{CCC}(G))}{|e(\mathcal{CCC}(G))|} \geq \frac{M_{1}(\mathcal{CCC}(G))}{|v(\mathcal{CCC}(G))|}$ with equality when $m=1$ and $2$.
\end{theorem}
\begin{proof}
     \textbf{Case 1.} $m$ is even.

     We have $\mathcal{CCC}(G)=K_{2m-2} \sqcup 2K_2$. Therefore, by Theorem \ref{thm1}, we get the expressions for $M_1(\mathcal{CCC}(G))$ and $M_2(\mathcal{CCC}(G))$.
     
     Note that $|v(\mathcal{CCC}(G))|=2m-2+4=2m+2$ and $|e(\mathcal{CCC}(G))|=\binom{2m-2}{2}+2=\frac{(2m-2)(2m-3)}{2}+2=(m-1)(2m-3)+2$. Therefore, 
\[
\frac{M_{2}(\mathcal{CCC}(G))}{|e(\mathcal{CCC}(G))|} = \frac{(m-1)(2m-3)^3+2}{(m-1)(2m-3)+2} \text{ and } \frac{M_{1}(\mathcal{CCC}(G))}{|v(\mathcal{CCC}(G))|} = \frac{(m-1)(2m-3)^2+2}{m+1}.
\]
Hence,
\begin{align*}
\frac{M_{2}(\mathcal{CCC}(G))}{|e(\mathcal{CCC}(G))|} - \frac{M_{1}(\mathcal{CCC}(G))}{|v(\mathcal{CCC}(G))|} 
&=\frac{8m^3(m-6)+8m(13m-12)+32}{(m+1)(m-1)(2m-3)+2(m+1)} 
:=\frac{f(m)}{g(m)}.
\end{align*}
Note that $g(m)=\frac{1}{2} \times |e(\mathcal{CCC}(G))| \times |v(\mathcal{CCC}(G))|>0$. We have $f(2)=0$ and $f(4)=288>0$. Also, since $13m-12>0$ for all $m \geq 6$ we have $\frac{f(m)}{g(m)}>0$.

\vspace{.2cm}

 \noindent \textbf{Case 2.} $m$ is odd.

 We have $\mathcal{CCC}(G)=K_{2m-1} \sqcup 2K_1$. Therefore, the result follows from the first case of Theorem \ref{CCC SD_8m}. 
 
  This completes the proof.
\end{proof}

\begin{theorem}
If $G=U_{(n,m)}$, where $m \geq 3$ and $n \geq 2$, then
     \begin{align*}
     M_1(\mathcal{CCC}(G)) = \begin{cases}
        \frac{1}{8}\left(16n(n-1)^2+(mn-2n)(mn-2n-2)^2\right),\\
         \hspace{6cm}\text{ when } m \text{ is even and } \frac{m}{2} \text{ is even}\\
        \frac{1}{8}\left(16n(2n-1)^2+(mn-2n)(mn-2n-2)^2\right),\\ \hspace{6cm} \text{ when } m \text{ is even and } \frac{m}{2} \text{ is odd}\\
        \frac{1}{8}\left(8n(n-1)^2+(mn-n)(mn-n-2)^2\right),
        \text{ when } m \text{ is odd}
    \end{cases}
    \end{align*}
    and
    \begin{align*}
    M_2(\mathcal{CCC}(G)) = \begin{cases}
        \frac{1}{32}\left(32n(n-1)^3+(mn-2n)(mn-2n-2)^3\right), \\
        \hspace{6cm} \text{ when } m \text{ is even and } \frac{m}{2} \text{ is even}\\
        \frac{1}{32}\left(32n(2n-1)^3+(mn-2n)(mn-2n-2)^3\right),\\ \hspace{6cm} \text{ when } m \text{ is even and } \frac{m}{2} \text{ is odd}\\
        \frac{1}{32}\left(16n(n-1)^3+(mn-n)(mn-n-2)^3\right),  \text{ when } m \text{ is odd.}
        \end{cases}
    \end{align*}
    Further, $\frac{M_{2}(\mathcal{CCC}(G))}{|e(\mathcal{CCC}(G))|} \geq \frac{M_{1}(\mathcal{CCC}(G))}{|v(\mathcal{CCC}(G))|}$ with equality when $m=3, 4$ and $6$.
\end{theorem}
\begin{proof}
    \textbf{Case 1.} $m$ is even and $\frac{m}{2}$ is even.

     We have $\mathcal{CCC}(G)=2K_n \sqcup K_{n(\frac{m}{2}-1)}$. Therefore, by Theorem \ref{thm1}, we get the expressions for $M_1(\mathcal{CCC}(G))$ and $M_2(\mathcal{CCC}(G))$.
     
     Note that  $|v(\mathcal{CCC}(G))|=2n+n(\frac{m}{2}-1)=\frac{1}{2}(mn+2n)$ and $|e(\mathcal{CCC}(G))|=2\binom{n}{2}+\binom{\frac{mn}{2}-n}{2}=\frac{1}{8}(8n(n-1)+(mn-2n)(mn-2n-2))$.  Therefore,
\[
\frac{M_{2}(\mathcal{CCC}(G))}{|e(\mathcal{CCC}(G))|} = \frac{32n(n-1)^3+(mn-2n)(mn-2n-2)^3}{4(8n(n-1)+(mn-2n)(mn-2n-2))} \text{ and }\]
\[
 \frac{M_{1}(\mathcal{CCC}(G))}{|v(\mathcal{CCC}(G))|} = \frac{16n(n-1)^2+(mn-2n)(mn-2n-2)^2}{4(mn+2n)}.
\]
Hence,
\begin{align*}
&\frac{M_{2}(\mathcal{CCC}(G))}{|e(\mathcal{CCC}(G))|} - \frac{M_{1}(\mathcal{CCC}(G))}{|v(\mathcal{CCC}(G))|} \\
&=\frac{1}{4}\times\frac{4n^4(m-2)(m-4)^2(mn-4)}{8n(n-1)(mn+2n)+(mn+2n)(mn-2n)(mn-2n-2)}
:=\frac{1}{4} \times \frac{f_1(m,n)}{g_1(m,n)}.
\end{align*}
Note that $g_1(m,n)=16 \times |e(\mathcal{CCC}(G))| \times |v(\mathcal{CCC}(G))|>0$. We have $f_1(4,n)=0$. Also, since $m-2>0$, $m-4>0$ and $mn-4>0$ for all $m \geq 6$ and $n \geq 2$, we have $\frac{1}{4} \times \frac{f_1(m,n)}{g_1(m,n)} \geq 0$. Hence the result follows.

\vspace{.3cm}

 \noindent \textbf{Case 2.} $m$ is even and $\frac{m}{2}$ is odd.

 We have $\mathcal{CCC}(G)=K_{2n} \sqcup K_{n(\frac{m}{2}-1)}$. Therefore, by Theorem \ref{thm1}, we get the expressions for $M_1(\mathcal{CCC}(G))$ and $M_2(\mathcal{CCC}(G))$.
 
 Note that $|v(\mathcal{CCC}(G))|=2n+n(\frac{m}{2}-1)=\frac{1}{2}(mn+2n)$ and $|e(\mathcal{CCC}(G))|=\binom{2n}{2}+\binom{\frac{mn}{2}-n}{2}=\frac{1}{8}(8n(2n-1)+(mn-2n)(mn-2n-2))$. Therefore,
 \[
 \frac{M_{2}(\mathcal{CCC}(G))}{|e(\mathcal{CCC}(G))|} = \frac{32n(2n-1)^3+(mn-2n)(mn-2n-2)^3}{4(8n(2n-1)+(mn-2n)(mn-2n-2))} \text{ and }
 \] 
 \[
 \frac{M_{1}(\mathcal{CCC}(G))}{|v(\mathcal{CCC}(G))|} = \frac{16n(2n-1)^2+(mn-2n)(mn-2n-2)^2}{4(mn+2n)}.
 \]
 Hence,
\begin{align*}
&\frac{M_{2}(\mathcal{CCC}(G))}{|e(\mathcal{CCC}(G))|} - \frac{M_{1}(\mathcal{CCC}(G))}{|v(\mathcal{CCC}(G))|}\\ 
&=\frac{1}{4}\times\frac{4n^4(m-2)(m-6)^2(mn+2n-4)}{8n(2n-1)(mn+2n)+(mn+2n)(mn-2n)(mn-2n-2)} 
:=\frac{1}{4}\times\frac{f_2(m,n)}{g_2(m,n)}.
\end{align*}
Note that $g_2(m,n)=16 \times |e(\mathcal{CCC}(G))| \times |v(\mathcal{CCC}(G))|>0$. We have $f_2(6,n)=0$. Also, since $m-2>0$, $m-6>0$ and $mn+2n-4>0$ for all $m \geq 10$ and $n \geq 2$, we have $\frac{1}{4}\times\frac{f_2(m,n)}{g_2(m,n)} \geq 0$. Hence the result follows.

\vspace{.2cm}

 \noindent \textbf{Case 3.} $m$ is odd.
                                         
  We have $\mathcal{CCC}(G)= K_n \sqcup K_{\frac{n(m-1)}{2}}$. Therefore, by Theorem \ref{thm1}, we get the expressions for $M_1(\mathcal{CCC}(G))$ and $M_2(\mathcal{CCC}(G))$.
  
  Note that $|v(\mathcal{CCC}(G))|=\frac{mn-n}{2}+n=\frac{1}{2}(mn+n)$ and $|e(\mathcal{CCC}(G))|=\binom{\frac{mn-n}{2}}{2}+\binom{n}{2}=\frac{1}{8}((mn-n)(mn-n-2)+4n(n-1)$. 
     Therefore,
\[
\frac{M_{2}(\mathcal{CCC}(G))}{|e(\mathcal{CCC}(G))|} = \frac{(mn-n)(mn-n-2)^3+16n(n-1)^3}{4((mn-n)(mn-n-2)+4n(n-1))} \text{ and }
\] 
\[
\frac{M_{1}(\mathcal{CCC}(G))}{|v(\mathcal{CCC}(G))|} = \frac{(mn-n)(mn-n-2)^2+8n(n-1)^2}{4(mn+n)}.
\]
Hence,
\begin{align*}
&\frac{M_{2}(\mathcal{CCC}(G))}{|e(\mathcal{CCC}(G))|} - \frac{M_{1}(\mathcal{CCC}(G))}{|v(\mathcal{CCC}(G))|} \\
&=\frac{1}{4}\times\frac{2m^3n^4(mn-6n-4)+16m^2n^5+4mn^4(14m-30)+6n^5(2m-3)+72n^4}{(mn+n)(mn-n)(mn-n-2)+4n(n-1)(mn+n)} \\
&:=\frac{1}{4}\times\frac{f_3(m,n)}{g_3(m,n)}.
\end{align*}
 Note that $g_3(m,n)=16 \times |e(\mathcal{CCC}(G))| \times |v(\mathcal{CCC}(G))|>0$. We have $f_3(3,n)=0$, $f_3(5,n)=192n^5-128n^4>0$, $f_3(7,n)=1536n^5-768n^4>0$ and $f_3(9,n)=5760n^5-2304n^4>0$ for all $n \geq 1$. Also, since $n(m-6)-4=mn-6n-4>0$, $14m-30>0$ and $2m-3>0$ for all $m \geq 11$ and $n \geq 1$, we have $\frac{1}{4}\times\frac{f_3(m,n)}{g_3(m,n)} \geq 0$. Hence, the result follows.
\end{proof}

We conclude this section with the following result.
\begin{theorem}
    If $G=G(p,m,n)$, where $p$ is a prime and $m, n \geq 1$, then
    \begin{align*}
        M_1(\mathcal{CCC}(G))&=2p^{3m+3n-3}-6p^{3m+3n-4}+6p^{3m+3n-5}-2p^{3m+3n-6}-4p^{2m+2n-2}+p^{3m+n-4} \\
        & \quad -4p^{2m+2n-4}+p^{3m+n}-4p^{3m+n-1}+6p^{3m+n-2}-4p^{3m+n-3}+8p^{2m+2n-3} \\
        & \quad -2p^{2m+n}+6p^{2m+n-1}-6p^{2m+n-2}+2p^{2m+n-3}+p^{m+n}-p^{m+n-2}
    \end{align*}
    and
    \begin{align*}
         M_2(\mathcal{CCC}(G))&=\frac{1}{2}\left(2p^{4m+4n-4}-8p^{4m+4n-5}+6p^{3m+3n-6}-8p^{4m+4n-7}+2p^{4m+4n-8}-p^{m+n}\right.\\
         & \quad \left.-6p^{3m+3n-3}+18p^{3m+3n-4}-18p^{3m+3n-5}+12p^{4m+4n-6}+6p^{2m+2n-2} \right.\\
         & \quad \left.+6p^{2m+2n-4}+p^{4m+n}-5p^{4m+n-1}+10p^{4m+n-2}-10p^{4m+n-3}+5p^{4m+n-4}\right.\\
         & \quad \left.+12p^{3m+n-1}-18p^{3m+n-2}+12p^{3m+n-3}-3p^{3m+n-4}+3p^{2m+n}-9p^{2m+n-1}\right. \\
         & \quad \left.-12p^{2m+2n-3}-3p^{3m+n}-p^{4m+n-5}+9p^{2m+n-2}-3p^{2m+n-3}+p^{m+n-2}\right).
    \end{align*}
Further, $\frac{M_{2}(\mathcal{CCC}(G))}{|e(\mathcal{CCC}(G))|} \geq \frac{M_{1}(\mathcal{CCC}(G))}{|v(\mathcal{CCC}(G))|}$ with equality when $n=1$.
\end{theorem}
\begin{proof}
    We have $\mathcal{CCC}(G)=K_{p^{m-1}(p^n-p^{n-1})} \sqcup K_{p^{n-1}(p^m-p^{m-1})} \sqcup (p^n-p^{n-1})K_{p^{m-n}(p^n-p^{n-1})}$ \\ $=2K_{(p^{m+n-1}-p^{m+n-2})} \sqcup (p^n-p^{n-1})K_{(p^m-p^{m-1})}$. Therefore,  by Theorem \ref{thm1}, we have
     \begin{align*}
         M_1(\mathcal{CCC}(G))=2(p^{m+n-1}-p^{m+n-2})&(p^{m+n-1}-p^{m+n-2}-1)^2 \\
         &+(p^n-p^{n-1})(p^m-p^{m-1})(p^m-p^{m-1}-1)^2 
     \end{align*}
     and
     \begin{align*}
         M_2(\mathcal{CCC}(G))=2(p^{m+n-1}-p^{m+n-2})&\times\frac{(p^{m+n-1}-p^{m+n-2}-1)^3}{2}\\
         &+(p^n-p^{n-1})(p^m-p^{m-1})\times\frac{(p^m-p^{m-1}-1)^3}{2}.
     \end{align*}
Hence, we get the required expressions for  $M_{1}(\mathcal{CCC}(G))$ and $M_{2}(\mathcal{CCC}(G))$ on simplification.

Note that $|v(\mathcal{CCC}(G))|=2(p^{m+n-1}-p^{m+n-2})+(p^n-p^{n-1})(p^m-p^{m-1})=p^{m+n}-p^{m+n-2}$ and $|e(\mathcal{CCC}(G))|=2 \binom{p^{m+n-1}-p^{m+n-2}}{2}+(p^n-p^{n-1})\binom{p^m-p^{m-1}}{2}$ $=\frac{1}{2}\left(2p^{2m+2n-2}-4p^{2m+2n-3}+\right. \\2p^{2m+2n-4}+p^{2m+n}-3p^{2m+n-1} \left.+3p^{2m+n-2}-p^{2m+n-3}-p^{m+n}+p^{m+n-2}\right)$.
We have

\noindent $\frac{M_{2}(\mathcal{CCC}(G))}{|e(\mathcal{CCC}(G))|} - \frac{M_{1}(\mathcal{CCC}(G))}{|v(\mathcal{CCC}(G))|}  =\frac{f(p,m,n)}{g(p,m,n)}$, where
     \begin{align*}
         f(p,m,n)&=2p^{5m+5n-4}-12p^{5m+5n-5}+30p^{5m+5n-6}-40p^{5m+5n-7}+30p^{5m+5n-8} \\
         & \quad-12p^{5m+5n-9}+2p^{5m+5n-10}-2p^{5m+4n-3}+12p^{5m+4n-4}-30p^{5m+4n-5} \\
         & \quad +40p^{5m+4n-6}-30p^{5m+4n-7}+12p^{5m+4n-8}-2p^{5m+4n-9}-4p^{4m+4n-3} \\
         & \quad +20p^{4m+4n-4}-40p^{4m+4n-5}+40p^{4m+4n-6}-20p^{4m+4n-7}+4p^{4m+4n-8}\\
         & \quad -2p^{5m+3n-2}+12p^{5m+3n-3}-30p^{5m+3n-4}+40p^{5m+3n-5}-30p^{5m+3n-6} \\
         & \quad +12p^{5m+3n-7}-2p^{5m+3n-8}+8p^{4m+3n-2}-40p^{4m+3n-3}+80p^{4m+3n-4} \\
         & \quad -80p^{4m+3n-5}+40p^{4m+3n-6}-8p^{4m+3n-7}+2p^{5m+2n-1}-12p^{5m+2n-2} \\
         & \quad +30p^{5m+2n-3}-40p^{5m+2n-4}+30p^{5m+2n-5}-12p^{5m+2n-6}+2p^{5m+2n-7} \\
         & \quad -4p^{4m+2n-1}+20p^{4m+2n-2}-40p^{4m+2n-3}+40p^{4m+2n-4}-20p^{4m+2n-5} \\
         & \quad \quad \quad \quad \quad \quad \quad \quad \quad \quad \quad \quad \quad \quad \quad \quad \quad \quad \quad \quad \quad \quad \quad \quad \quad +4p^{4m+2n-6}.
         \end{align*}
     and \begin{align*}
         g(p,m,n)&=(p^{m+n}-p^{m+n-2})(2p^{2m+2n-2}-4p^{2m+2n-3}+2p^{2m+2n-4}+p^{2m+n}-3p^{2m+n-1} \\
         & \quad \quad \quad \quad \quad \quad \quad \quad \quad \quad \quad \quad \quad \quad \quad  +3p^{2m+n-2}-p^{2m+n-3}-p^{m+n}+p^{m+n-2}).
     \end{align*}
     On factorization, we have $f(p,m,n)=2p^{4m+2n-10}(p-1)^5  (p - p^n)^2(p^m(p^{n+1}-p^n-p)+p^2(p^m-2))$. Clearly, $f(p,m,1)=0$. For all prime $p \geq 2, m \geq 1$ and $n > 1$, we have $p^{n+1}-p^n-p>0$ and $p^m-2 \geq 0$. As such, $f(p,m,n) > 0$ and $g(p,m,n)=2 \times |e(\mathcal{CCC}(G))| \times |v(\mathcal{CCC}(G))| > 0$. Thus, $\frac{f(p,m,n)}{g(p,m,n)} \geq 0$ and the result follows.
\end{proof}
\subsection{Zagreb indices of CCC-graphs of groups with a given central quotient}
In this section, we consider the CCC-graphs of the groups  $G$ when $G/Z(G)$ is isomorphic to   $D_{2m}$, $\mathbb{Z}_p \times \mathbb{Z}_p$, Frobenious group of order $pq$ and $p^2q$ or any  group of order $p^3$ for any two primes $p$ and $q$. 
The structures of CCC-graphs of groups when its central quotient is isomorphic to a dihedral group have been obtained in \cite[Theorem 1.2]{MAS-2020}. Using which we get the following result.
\begin{theorem}\label{cen-iso-D}
    Let $G$ be a finite group such that $\frac{G}{Z(G)} \cong D_{2m}$, where $m \geq 3.$ Then  
    \begin{align*}
             M_1(\mathcal{CCC}(G))= \begin{cases}
                 \frac{1}{8}\left((mx-x)(mx-x-2)^2+2x(x-2)^2\right), & \text{ when } m \text{ is even} \\
                 \frac{1}{8}\left((mx-x)(mx-x-2)^2+8x(x-1)^2\right), & \text{ when } m \text{ is odd}
             \end{cases}
    \end{align*}
    and
    \begin{align*}
         M_2(\mathcal{CCC}(G))= \begin{cases}
                 \frac{1}{32}\left((mx-x)(mx-x-2)^3+2x(x-2)^3\right), & \text{ when } m \text{ is even} \\
                 \frac{1}{32}\left((mx-x)(mx-x-2)^3+16x(x-1)^3\right), & \text{ when } m \text{ is odd,}
             \end{cases}
    \end{align*}
    where $x=|Z(G)|$. Further, $\frac{M_{2}(\mathcal{CCC}(G))}{|e(\mathcal{CCC}(G))|} > \frac{M_{1}(\mathcal{CCC}(G))}{|v(\mathcal{CCC}(G))|}$.
\end{theorem}
\begin{proof}
    \textbf{Case 1.} $m$ is even.

     We have $\mathcal{CCC}(G)=K_{\frac{x(m-1)}{2}} \sqcup 2K_{\frac{x}{2}}$. Therefore, by Theorem \ref{thm1}, we get the expressions for $M_1(\mathcal{CCC}(G))$ and $M_2(\mathcal{CCC}(G))$.
     
     Note that $|v(\mathcal{CCC}(G))|=\frac{x(m-1)}{2}+2 \times \frac{x}{2}=\frac{x(m+1)}{2}$ and $|e(\mathcal{CCC}(G))|=\binom{\frac{mx-x}{2}}{2}+2\binom{\frac{x}{2}}{2}=\frac{1}{8}((mx-x)(mx-x-2)+2x(x-2))$.
Therefore,
\[
\frac{M_{2}(\mathcal{CCC}(G))}{|e(\mathcal{CCC}(G))|} = \frac{(mx-x)(mx-x-2)^3+2x(x-2)^3}{4((mx-x)(mx-x-2)+2x(x-2))} \text{ and }
\]
\[
\frac{M_{1}(\mathcal{CCC}(G))}{|v(\mathcal{CCC}(G))|} 
= \frac{(mx-x)(mx-x-2)^2+2x(x-2)^2}{4x(m+1)}.
\]
Hence, 
\begin{align*}
\frac{M_{2}(\mathcal{CCC}(G))}{|e(\mathcal{CCC}(G))|} - \frac{M_{1}(\mathcal{CCC}(G))}{|v(\mathcal{CCC}(G))|} 
&=\frac{1}{4}\times\frac{x^2(m-2)^2(mx-4)}{mx(mx+x)(mx-2x-2)}
>0 
\end{align*}
for all $m \geq 4$.

\vspace{.2cm}

 \noindent \textbf{Case 2.} $m$ is odd.
                                         
  We have $\mathcal{CCC}(G)=K_{\frac{x(m-1)}{2}} \sqcup K_x$. Therefore, by Theorem \ref{thm1}, we get the expressions for $M_1(\mathcal{CCC}(G))$ and $M_2(\mathcal{CCC}(G))$.
  
  Note that $|v(\mathcal{CCC}(G))|=\frac{mx-x}{2}+x=\frac{1}{2}(mx+x)$ and $|e(\mathcal{CCC}(G))|=\binom{\frac{mx-x}{2}}{2}+\binom{x}{2}=\frac{1}{8}((mx-x)(mx-x-2)+4x(x-1))$. Therefore,
  \[
  \frac{M_{2}(\mathcal{CCC}(G))}{|e(\mathcal{CCC}(G))|} = \frac{(mx-x)(mx-x-2)^3+16x(x-1)^3}{4((mx-x)(mx-x-2)+4x(x-1))} \text{ and }
  \]
  \[
  \frac{M_{1}(\mathcal{CCC}(G))}{|v(\mathcal{CCC}(G))|} 
  = \frac{(mx-x)(mx-x-2)^2+8x(x-1)^2}{4(mx+x)}.
  \]
  Hence,
\begin{align*}
&\frac{M_{2}(\mathcal{CCC}(G))}{|e(\mathcal{CCC}(G))|} - \frac{M_{1}(\mathcal{CCC}(G))}{|v(\mathcal{CCC}(G))|}\\ 
&=\frac{1}{4}\times\frac{2m^3x^4(mx-6x-4)+16m^2x^5+4mx^4(14m-30)+6x^5(2m-3)+72x^4}{(mx+x)(mx-x)(mx-x-2)+4x(x-1)(mx+x)} \\
&:=\frac{1}{4}\frac{f(m,x)}{g(m,x)}.
\end{align*}
 Note that $g(m,x)=16 \times |e(\mathcal{CCC}(G))| \times |v(\mathcal{CCC}(G))|>0$. We have $f(3,x)=0$, $f(5,x)=192x^5-128x^4>0$, $f(7,x)=1536x^5-768x^4>0$ and $f(9,x)=5760x^5-2304x^4>0$ for all $x \geq 1$. Also, since $x(m-6)-4=mx-6x-4>0$, $14m-30>0$ and $2m-3>0$ for all $m \geq 11$ and $x \geq 1$, we have $\frac{1}{4}\frac{f(m,x)}{g(m,x)} \geq 0$. Hence, the result follows.
\end{proof}

The structures of CCC-graphs of groups when its central quotient is isomorphic to $\mathbb{Z}_p \times \mathbb{Z}_p$ (for any prime $p$) have been obtained in \cite[Theotem 3.1]{MAS-ARA-2020}. Therefore,  we get the following result.

\begin{theorem}\label{G/Z=Z_p x Z_p}
    Let $G$ be a finite  group   and $\frac{G}{Z(G)} \cong \mathbb{Z}_p \times \mathbb{Z}_p$, where $p$ is a prime. Then
    \[
    M_1(\mathcal{CCC}(G))= \frac{x(p^2-1)}{p^3}\left(p^2x^2-2px^2+x^2-2p^2x+2px+p^2\right) \text{ \quad and }
    \]
\begin{align*}
M_2&(\mathcal{CCC}(G))\\
&=
\frac{x(p^2-1)}{2p^4}\left(p^3x^3-3p^2x^3+3px^3-x^3-3p^3x^2+6p^2x^2-3px^2\right.  \left. +3p^3x-3p^2x-p^3\right),
\end{align*}  
    where $x=|Z(G)|$. Further, $\frac{M_{2}(\mathcal{CCC}(G))}{|e(\mathcal{CCC}(G))|} = \frac{M_{1}(\mathcal{CCC}(G))}{|v(\mathcal{CCC}(G))|}$.
\end{theorem}
\begin{proof}
    We have $\mathcal{CCC}(G) = (p+1)K_{\frac{x(p-1)}{p}}$. Therefore, by Theorem \ref{thm1}, we get the expressions for $M_1(\mathcal{CCC}(G))$ and $M_2(\mathcal{CCC}(G))$.
    
    Note that $|v(\mathcal{CCC}(G))|=(p+1) \times \frac{x(p-1)}{p}=\frac{x(p^2-1)}{p}$ and $|e(\mathcal{CCC}(G))|=(p+1)\binom{\frac{x(p-1)}{p}}{2}=\frac{x(p^2-1)}{2p^2}(px-x-p)$. 
    We have
    \begin{align*}
       \frac{M_{2}(\mathcal{CCC}(G))}{|e(\mathcal{CCC}(G))|}
       &=\frac{p^2x^2-2px^2+x^2-2p^2x+2px+p^2}{p^2} = \frac{M_{1}(\mathcal{CCC}(G))}{|v(\mathcal{CCC}(G))|}.
      \end{align*}
     Hence, the result follows.
\end{proof}
\begin{cor}\label{order G/Z=p^2}
    If $G$ is a group of order $p^n$ and $|Z(G)|=p^{n-2}$, where $p$ is a prime and $n \geq 3$, then
    \begin{align*}
    &M_1(\mathcal{CCC}(G))\\
    &=p^{3n-5}-2p^{3n-6}+2p^{3n-8}-p^{3n-9}-2p^{2n-3}+2p^{2n-4}+2p^{2n-5}-2p^{2n-6} +p^{n-1}-p^{n-3}  
    \end{align*} 
and
    \begin{align*}
        M_2(\mathcal{CCC}(G))&=p^{4n-7}-3p^{4n-8}+2p^{4n-9}+2p^{4n-10}-3p^{4n-11}+p^{4n-12}-3p^{3n-5}+6p^{3n-6} \\
        &\quad \quad -6p^{3n-8}+3p^{3n-9}+3p^{2n-3}-3p^{2n-4}-3p^{2n-5}+3p^{2n-6}-p^{n-1}+p^{n-3}.
    \end{align*}
    Further,  $\frac{M_{2}(\mathcal{CCC}(G))}{|e(\mathcal{CCC}(G))|} = \frac{M_{1}(\mathcal{CCC}(G))}{|v(\mathcal{CCC}(G))|}$.
\end{cor}
\begin{proof}
    Clearly, $\frac{G}{Z(G)}\cong \mathbb{Z}_p \times \mathbb{Z}_p$ and $\frac{(p-1)|Z(G)|}{p}=(p-1)p^{n-3}$. Therefore, the result follows from Theorem \ref{G/Z=Z_p x Z_p}.
\end{proof}
Putting $n =3$ in Corollary \ref{order G/Z=p^2} we get the following result.
\begin{cor}
If $G$ is a non-abelian group of order $p^3$ then 
\[ 
	M_1(\mathcal{CCC}(G))
	=p^4 - 4p^{3} + 3p^2 + 4p - 4, \,   
	M_2(\mathcal{CCC}(G))= p^5 - 6p^4 + 2p^3 + 9p^3 - 2p^2
	-12p  +8
\]
and  $\frac{M_{2}(\mathcal{CCC}(G))}{|e(\mathcal{CCC}(G))|} = \frac{M_{1}(\mathcal{CCC}(G))}{|v(\mathcal{CCC}(G))|}$.
\end{cor}

Now we consider finite groups whose central quotient is isomorphic to certain Frobenious groups. The CCC-graphs of these groups  have been obtained in \cite{MR-ZF-2024}.
\begin{theorem}
Let $p$ and $q$ be two primes.    Let $G$ be a finite group such that $\frac{G}{Z(G)} \cong \mathbb{Z}_q \rtimes_f \mathbb{Z}_p$,  the Frobenious group of order $pq$ with kernel $\mathbb{Z}_q$ and complement $\mathbb{Z}_p$. Then
    \begin{align*}
        M_1(\mathcal{CCC}(G))&=\frac{x}{p^3}\left(q^3x^2-3q^2x^2+3qx^2+p^6x^2-3p^5x^2+3p^4x^2-p^3x^2-x^2-2pq^2x \right.\\
        & \quad \quad \quad \left. +4pqx-2p^5x+4p^4x-2p^3x-2px+p^2q+p^4-p^3-p^2 \right) \text{ \quad and }
    \end{align*}
    \begin{align*}
        M_2(\mathcal{CCC}(G))&=\frac{x}{2p^4}\left(q^4x^3-4q^3x^3+6q^2x^3-4px^3-p^3q-4p^7x^3+6p^6x^3-4p^5x^3+3px^2\right.\\
        & \quad \left.+3p^2q^2x-3pq^3x^2+9pq^2x^2-9pqx^2-3p^7x^2+9p^6x^2-9p^5x^2+3p^4x^2 \right.\\
        & \quad \left.+p^4x^3-6p^2qx+3p^6x-6p^5x+3p^4x+3p^2x+x^3+p^8x^3-p^5+p^4+p^3 \right)
    \end{align*}
    where $x=|Z(G)|$. Further,  $\frac{M_{2}(\mathcal{CCC}(G))}{|e(\mathcal{CCC}(G))|} \geq \frac{M_{1}(\mathcal{CCC}(G))}{|v(\mathcal{CCC}(G))|}$.
\end{theorem}
\begin{proof}
    We have $\mathcal{CCC}(G)=K_{\frac{x(q-1)}{p}} \sqcup K_{x(p-1)}$, where $x=|Z(G)|$.Therefore, by Theorem \ref{thm1}, we get the expressions for $M_1(\mathcal{CCC}(G))$ and $M_2(\mathcal{CCC}(G))$.
    
    Note that $|v(\mathcal{CCC}(G))|=\frac{x(q-1)}{p}+x(p-1)=\frac{x}{p}(p^2+q-p-1)$ and $|e(\mathcal{CCC}(G))|=\binom{\frac{x(q-1)}{p}}{2}+\binom{x(p-1)}{2}=\frac{x}{2p^2}(p^2x+q^2x+p^4x+x-2qx-pq-2p^3x-p^3+p^2+p)$. 
    It can be seen that
    \begin{align*}
    &\frac{M_{2}(\mathcal{CCC}(G))}{|e(\mathcal{CCC}(G))|} - \frac{M_{1}(\mathcal{CCC}(G))}{|v(\mathcal{CCC}(G))|} 
    =\frac{f(p,q,x)}{g(p,q,x)},
    \end{align*}
    where 
    \begin{align*}
       f(p,q,x)&=p^2q^4x^3-pq^4x^3-p^4q^3x^3+2p^3q^3x^3-5p^2q^3x^3+4pq^3x^3-p^6q^2x^3+3p^5q^2x^3 \\
       & \quad -5p^3q^2x^3+9p^2q^2x^3-6pq^2x^3+p^8qx^3-4p^7qx^3+8p^6qx^3-10p^5qx^3-p^3x^3 \\
       & \quad +4p^4qx^3+4p^3qx^3-7p^2qx^3+4pqx^3-p^8x^3+4p^7x^3-7p^6x^3+7p^5x^3-px^3 \\
       & \quad +2p^2x^3-2p^3q^3x^2+2p^2q^3x^2+4p^5q^2x^2-8p^4q^2x^2-6p^2q^2x^2-2p^7qx^2 \\
       & \quad +6p^6qx^2-14p^5qx^2+18p^4qx^2-14p^3qx^2+6p^2qx^2+2p^7x^2-6p^6x^2 \\
       & \quad \quad \quad \quad \quad \quad \quad \quad \quad  \quad \quad +10p^5x^2-3p^4x^3-10p^4x^2+6p^3x^2-2p^2x^2
    \end{align*}
and $g(p,q,x) = p^2(p^2+q-p-1)(p^2x+q^2x+p^4x+x-2qx-pq-2p^3x-p^3+p^2+p)$.
    On factorization, we have $f(p,q,x)=px^2(p-1)(q-1)(p+q-p^2-1)^2(p^2x+qx-x-2p-px)$. Suppose that $p^2x+qx-x-2p-px=p(x(p-1)-2)+x(q-1):=S_1 + S_2$, where $S_1=p(x(p-1)-2)$ and $S_2=x(q-1) >0$ for all prime $q$ and $x \geq 1$. Clearly, $S_1 \geq 0$ for all prime $p \geq 3$ and $x \geq 1$. Now, $f(2,q,1)=2(q-1)(q-3)^2(q-3) \geq 0$ for all $q \geq 3$. Also, $g(p,q,x)=\frac{2p^5}{x^2} \times |e(\mathcal{CCC}(G))| \times |v(\mathcal{CCC}(G))|>0$. Therefore, $\frac{f(p,q,x)}{g(p,q,x)} \geq 0$ and hence the result follows.
\end{proof}

\begin{theorem}
    Let $G$ be a finite group such that $\frac{G}{Z(G)}$ is isomorphic to a Frobenious group of order $p^2q$, where $p$ and $q$ are primes.
    Then $\mathcal{CCC}(G)$ satisfies Hansen-Vuki{\v{c}}evi{\'c} conjecture if $\frac{G}{Z(G)} \cong A_4$ or $p<q$ and $|\frac{G}{Z(G)}|\neq 12$ or  $p>q$.
\end{theorem} 
\begin{proof}
   We have the following three cases.
   
    \noindent \textbf{Case 1.} $\frac{G}{Z(G)} \cong A_4$.
    
     We have $\mathcal{CCC}(G)=K_{2x} \sqcup K_x$ or $K_{2x} \sqcup K_{\frac{x}{2}}$, where $x=|Z(G)|$. If $\mathcal{CCC}(G)=K_{2x} \sqcup K_x$ then, by Theorem \ref{thm1}, we have
        \begin{align*}
            M_1(\mathcal{CCC}(G))&=2x(2x-1)^2+x(x-1)^2=9x^3-10x^2+3x \text{ \quad and}
        \end{align*}
        \begin{align*}
            M_2(\mathcal{CCC}(G))&=2x\frac{(2x-1)^3}{2}+x\frac{(x-1)^3}{2}=\frac{1}{2}\left(17x^4-27x^3+15x^2-3x\right).
        \end{align*}
        Also, $|v(\mathcal{CCC}(G))|=2x+x=3x$ and $|e(\mathcal{CCC}(G))|=\binom{2x}{2}+\binom{x}{2}=\frac{1}{2}(5x^2-3x)$. Therefore, for all $x \geq 1$, we have
        \begin{align*}
            \frac{M_{2}(\mathcal{CCC}(G))}{|e(\mathcal{CCC}(G))|} - \frac{M_{1}(\mathcal{CCC}(G))}{|v(\mathcal{CCC}(G))|}=\frac{6x^5-4x^4}{3x(5x^2-3x)}>0.
        \end{align*}
        If $\mathcal{CCC}(G)=K_{2x} \sqcup K_{\frac{x}{2}}$ then, by Theorem \ref{thm1}, we have
        \begin{align*}
            M_1(\mathcal{CCC}(G))&=2x(2x-1)^2+\frac{x}{2}\left(\frac{x}{2}-1\right)^2=\frac{1}{8}\left(65x^3-68x^2+20x\right) \text{ \quad and}
        \end{align*}
        \begin{align*}
            M_2(\mathcal{CCC}(G))&=2x\frac{(2x-1)^3}{2}+\frac{x}{2}\frac{(\frac{x}{2}-1)^3}{2}=\frac{1}{32}\left(257x^4-390x^3+204x^2-40x\right).
        \end{align*}
        Also, $|v(\mathcal{CCC}(G))|=2x+\frac{x}{2}=\frac{5x}{2}$ and $|e(\mathcal{CCC}(G))|=\binom{2x}{2}+\binom{\frac{x}{2}}{2}=\frac{1}{8}(17x^2-10x)$. Therefore, for all $x \geq 1$, we have
        \begin{align*}
            \frac{M_{2}(\mathcal{CCC}(G))}{|e(\mathcal{CCC}(G))|} - \frac{M_{1}(\mathcal{CCC}(G))}{|v(\mathcal{CCC}(G))|}=\frac{180x^5-144x^4}{20x(17x^2-10x)}>0.
        \end{align*}

\vspace{.2cm}
    
\noindent        \textbf{Case 2.} $p<q$ and $|\frac{G}{Z(G)}|\neq 12$.
        
         We have $\mathcal{CCC}(G)=K_{\frac{x(q-1)}{p^2}}\sqcup K_{x(p^2-1)}$, where $x=|Z(G)|$. As such, $|v(\mathcal{CCC}(G))|=\frac{x(q-1)}{p^2}+x(p^2-1)=\frac{1}{p^2}(p^4x+qx-p^2x-x)$ and $|e(\mathcal{CCC}(G))|=\binom{\frac{x(q-1)}{p^2}}{2}+\binom{x(p^2-1)}{2}=\frac{1}{2p^4}(
        q^2x^2-2qx^2+p^8x^2-2p^6x^2+p^4x^2+x^2-p^2qx-p^6x+p^4x+p^2x)$. By Theorem \ref{thm1}, we have
        \begin{align*}
            M_1(\mathcal{CCC}(G))&=\frac{x(q-1)}{p^2}\left(\frac{x(q-1)}{p^2}-1\right)^2+x(p^2-1)\left(x(p^2-1)-1\right)^2 \\
            &=\frac{1}{p^6}\left(q^3x^3-3q^2x^3+3q x^3+p^{12}x^3-3p^{10}x^3+3p^8x^3-p^6x^3-2p^2q^2x^2 \right. \\
            & \quad \quad \quad \quad \quad \quad \quad \quad \left. -x^3+4p^2qx^2+4p^8x^2-2p^{10}x^2-2p^6x^2+p^4qx \right. \\
            & \quad \quad \quad \quad \quad \quad \quad \quad \quad \quad \quad \quad \left.-2p^2x^2+p^8x-p^4x-p^6x\right) \text{ \quad and }
            \end{align*} 
            \begin{align*}
            M_2(\mathcal{CCC}(G))&=\frac{x(q-1)}{p^2}\frac{\left(\frac{x(q-1)}{p^2}-1\right)^3}{2}+x(p^2-1)\frac{\left(x(p^2-1)-1\right)^3}{2} \\
            &=\frac{1}{2p^8}\left(q^4x^4-4q^3x^4+6q^2x^4-4qx^4+p^{16}x^4-4p^{14}x^4-4p^{10}x^4+6p^{12}x^4 \right. \\
            & \quad \quad \quad \quad \left.+x^4+p^8x^4-3p^2q^3x^3+9p^2q^2x^3-9p^2qx^3+9p^{12}x^3-3p^{14}x^3 \right.\\
            & \quad \quad \quad \quad \quad \left.+3p^8x^3-9p^{10}x^3+3p^4q^2x^2+3p^2x^3-6p^4qx^2+3p^{12}x^2\right. \\
            & \quad \quad \quad \quad \quad \left.-6p^{10}x^2+3p^4x^2 +3p^8x^2-p^{10}x-p^6qx+p^6x+p^8x\right).
            \end{align*}
            Therefore, $\frac{M_{2}(\mathcal{CCC}(G))}{|e(\mathcal{CCC}(G))|} - \frac{M_{1}(\mathcal{CCC}(G))}{|v(\mathcal{CCC}(G))|}=\frac{f_1(p,q,x)}{g_1(p,q,x)}$, where $f_1(p,q,x)=p^2x^4(p-1)(p+1)(p^4-p^2-q+1)^2(q-1)(p^2(p^2x-x-2)+x(q-1))$ after factorization and $g_1(p,q,x)=2p^{10} \times |e(\mathcal{CCC}(G))| \times |v(\mathcal{CCC}(G))| >0$. 
            For all prime $p \geq 2$ and $x \geq 1$, we have $x(p^2-1) > 2$. As such, $f_1(p,q,x) >0$ and therefore, $\frac{f_1(p,q,x)}{g_1(p,q,x)} >0$.
            
            \vspace{.2cm}

            \noindent \textbf{Case 3.} $p > q$.
            
             We have $\mathcal{CCC}(G)=K_{x(q-1)}\sqcup K_{\frac{x(p^2-1)}{q}}$ or $K_{x(q-1)}\sqcup (p+1)K_{\frac{x(p-1)}{pq}}$. If $\mathcal{CCC}(G)=K_{x(q-1)}\sqcup K_{\frac{x(p^2-1)}{q}}$ then, by Theorem \ref{thm1}, we have
            \begin{align*}
                M_{1}&(\mathcal{CCC}(G))=x(q-1)\left(x(q-1)-1\right)^2+\frac{x(p^2-1)}{q}\left(\frac{x(p^2-1)}{q}-1\right)^2 \\
                &\qquad=\frac{1}{q^3}\left(q^6x^3-3q^5x^3+3q^4x^3-q^3x^3+p^6x^3-3p^4x^3+3p^2x^3-2q^5x^2-x^3 \right. \\
                &\qquad \quad \quad  \quad \left.+4q^4x^2-2q^3x^2-2p^4qx^2+4p^2qx^2-2qx^2+q^4x -q^3x+p^2q^2x-q^2x\right) 
            \end{align*}
        and
            \begin{align*}
                M_{2}(\mathcal{CCC}(G))&=x(q-1)\frac{\left(x(q-1)-1\right)^3}{2}+\frac{x(p^2-1)}{q}\frac{\left(\frac{x(p^2-1)}{q}-1\right)^3}{2} \\
                &=\frac{1}{2q^4}\left(q^8x^4-4q^7x^4+6q^6x^4-4q^5x^4+p^8x^4+q^4x^4-4p^6x^4+6p^4x^4 \right. \\
                & \quad \quad \quad \left.-4p^2x^4+x^4-3q^7x^3+9q^6x^3-9q^5x^3+3q^4x^3-3p^6qx^3+9p^4qx^3 \right. \\
                & \quad \quad \quad \left.-9p^2qx^3+3q^6x^2+3qx^3-6q^5x^2+3p^4q^2x^2+3q^4x^2-6p^2q^2x^2 \right. \\
                &\quad \quad \quad \quad \quad \quad \quad \quad \quad \quad \quad \quad \quad \left.+3q^2x^2-q^5x+q^4x+q^3x-p^2q^3x\right).
            \end{align*}
            Also, $|v(\mathcal{CCC}(G))|=x(q-1)+\frac{x(p^2-1)}{q}=\frac{1}{q}(q^2x-qx+p^2x-x)$ and $|e(\mathcal{CCC}(G))|=\binom{x(q-1)}{2}+\binom{\frac{x(p^2-1)}{q}}{2}=\frac{1}{2q^2}(q^4x^2-2q^3x^2+q^2x^2+p^4x^2-2p^2x^2-q^3x+x^2-p^2qx+q^2x+qx)$. Therefore, $\frac{M_{2}(\mathcal{CCC}(G))}{|e(\mathcal{CCC}(G))|} - \frac{M_{1}(\mathcal{CCC}(G))}{|v(\mathcal{CCC}(G))|}=\frac{f_2(p,q,x)}{g_2(p,q,x)}$, where $f_2(p,q,x)=qx^4(p-1)(p+1)(q-1)(q^2-p^2-q+1)^2(x(p^2-q-1)+q(qx-2))$ after factorization and $g_2(p,q,x)=2q^5 \times |e(\mathcal{CCC}(G))| \times |v(\mathcal{CCC}(G))| >0$. Since $p>q$, so $p^2-q>1$ and $qx>2$ for all prime $q \geq 2$ and $x \geq 1$. As such, $f_2(p,q,x)>0$ and therefore, $\frac{f_2(p,q,x)}{g_2(p,q,x)}>0$.

             If $\mathcal{CCC}(G)=K_{x(q-1)}\sqcup (p+1)K_{\frac{x(p-1)}{pq}}$ then, by Theorem \ref{thm1}, we have
             \begin{align*}
                 M_{1}(\mathcal{CCC}(G))&=x(q-1)\left(x(q-1)-1\right)^2+(p+1)\frac{x(p-1)}{pq}\left(\frac{x(p-1)}{pq}-1\right)^2 \\
                 &=\frac{1}{p^3q^3}\left(p^3q^6x^3-3p^3q^5x^3+3p^3q^4x^3-p^3q^3x^3+p^4x^3-2p^3x^3+2px^3\right. \\
                 & \quad \quad \quad \left. -2p^3q^5x^2-x^3+4p^3q^4x^2-2p^4qx^2-2p^3q^3x^2+2p^2qx^2\right. \\
                 & \quad \quad \quad \left.+2p^3qx^2+p^3q^4x-2pqx^2+p^4q^2x-p^3q^3x-p^2q^2x\right) 
             \end{align*}
  and
             \begin{align*}
                 M_{2}(\mathcal{CCC}(G))&=x(q-1)\frac{\left(x(q-1)-1\right)^3}{2}+(p+1)\frac{x(p-1)}{pq}\frac{\left(\frac{x(p-1)}{pq}-1\right)^3}{2} \\
                 &=\frac{1}{2p^4q^4}\left(p^4q^8x^4-4p^4q^7x^4+6p^4q^6x^4-4p^4q^5x^4+p^4q^4x^4+p^5x^4-3p^4x^4 \right. \\
                 & \quad \quad \quad \quad \left.+2p^2x^4+2p^3x^4-3px^4+x^4-3p^4q^7x^3+9p^4q^6x^3-9p^4q^5x^3\right.\\
                 &\quad \quad \quad \quad \left.+3p^4q^4x^3-3p^5qx^3+6p^4qx^3-6p^2qx^3+3p^4q^6x^2+3pqx^3 \right.\\
                 &\quad \quad \quad \quad \left.- 6p^4q^5x^2+3p^5q^2x^2+3p^4q^4x^2-3p^3q^2x^2-3p^4q^2x^2+3p^2q^2x^2\right. \\
                 & \quad \quad \quad \quad \quad \quad \quad \quad \quad \quad \quad \quad \quad \quad \quad \left.-p^4q^5x+p^4q^4x+p^3q^3x-p^5q^3x\right).
             \end{align*}
Also, $|v(\mathcal{CCC}(G))|=x(q-1)+(p+1)\frac{x(p-1)}{pq}=\frac{1}{pq}(pq^2x-pqx+p^2x-x)$ and $|e(\mathcal{CCC}(G))|=\binom{x(q-1)}{2}+(p+1)\binom{\frac{x(p-1)}{pq}}{2}=\frac{1}{2p^2q^2}(p^2q^4x^2-2p^2q^3x^2+p^2q^2x^2+p^3x^2-px^2-p^2x^2+x^2-p^2q^3x+p^2q^2x+pqx-p^3qx)$. Therefore, $\frac{M_{2}(\mathcal{CCC}(G))}{|e(\mathcal{CCC}(G))|} - \frac{M_{1}(\mathcal{CCC}(G))}{|v(\mathcal{CCC}(G))|}=\frac{f_3(p,q,x)}{g_3(p,q,x)}$, where $f_3(p,q,x)=pqx^4(p-1) (p+1)(q-1)(pq^2-pq-p+1)^2(pq(qx-x-2)+x(p-1))$ after factorization and $g_3(p,q,x)=2p^5q^5 \times |e(\mathcal{CCC}(G))| \times |v(\mathcal{CCC}(G))| >0$. Since $x\geq p > q$, so $x(q-1)>2$. As such, $f_3(p,q,x)>0$ and therefore, $\frac{f_3(p,q,x)}{g_3(p,q,x)}>0$. 

In all the cases, we have $\frac{M_{2}(\mathcal{CCC}(G))}{|e(\mathcal{CCC}(G))|} > \frac{M_{1}(\mathcal{CCC}(G))}{|v(\mathcal{CCC}(G))|}$ and hence the result follows.
\end{proof}
\begin{theorem}\label{G/Z abelian}
	Let $G$ be a non-abelian group  such that $\frac{G}{Z(G)}$ is abelian and $|\frac{G}{Z(G)}|=p^3$, where   $p$ is a prime. Then $\mathcal{CCC}(G)$ satisfies Hansen-Vuki{\v{c}}evi{\'c} conjecture. 
\end{theorem}
\begin{proof}
	We have $\mathcal{CCC}(G)=K_m \sqcup p^2K_n$ or $(p^2+p+1)K_n$, where $m=\frac{x(p^2-1)}{p}, n=\frac{x(p-1)}{p^2}$ and $x=|Z(G)|$. If $\mathcal{CCC}(G)=K_m \sqcup p^2K_n$ then, by Theorem \ref{thm1}, we have

	\begin{align*}
		M_1(\mathcal{CCC}(G))&= \frac{x(p^2-1)}{p}\left(\frac{x(p^2-1)}{p}-1\right)^2+p^2 \times \frac{x(p-1)}{p^2}\left(\frac{x(p-1)}{p^2}-1\right)^2 \\
		&=\frac{x(p-1)}{p^4}\left(p^6x^2-2p^4x^2+2p^2x^2-2p^5x+2p^4+p^5x^2-2p^3x^2-px^2 \right.\\
		& \quad \quad \quad \quad \quad \quad \quad \quad \quad \quad \quad \quad \quad \quad \quad \quad \left. -2p^4x+4p^2x+p^3+x^2\right) \quad \text{ and }
	\end{align*}
	\begin{align*}
		M_2(\mathcal{CCC}(G))&= \frac{x(p^2-1)}{p}\frac{\left(\frac{x(p^2-1)}{p}-1\right)^3}{2}+p^2 \times \frac{x(p-1)}{p^2}\frac{\left(\frac{x(p-1)}{p^2}-1\right)^3}{2} \\
		&=\frac{x(p-1)}{2p^6}\left(p^9x^3+p^8x^3-3p^8x^2-3p^7x^3-3p^7x^2-3p^6x^3+6p^6x^2\right.\\
		&\quad \quad \quad \quad \quad \quad \left. +3p^7x +3p^6x+3p^5x^3+6p^5x^2+3p^4x^3-6p^4x^2-6p^4x \right. \\
		& \quad \quad \quad \quad \quad \quad \quad \quad \quad \left.+3p^3x^2-4p^2x^3-p^5-2p^6-x^3-3p^2x^2+3px^3\right).
	\end{align*}
	Also, $|v(\mathcal{CCC}(G))|=\frac{x(p^2-1)}{p}+p^2 \times \frac{x(p-1)}{p^2}=x(p-1)(p+2)$ and $|e(\mathcal{CCC}(G))|=\binom{\frac{x(p^2-1)}{p}}{2}+p^2\binom{\frac{x(p-1)}{p^2}}{2}=\frac{x(p-1)}{2p^2}(p^3x+p^2x-2p^2-2x-p)$. Therefore,
	\begin{align*}
		&\frac{M_{2}(\mathcal{CCC}(G))}{|e(\mathcal{CCC}(G))|} -\frac{M_{1}(\mathcal{CCC}(G))}{|v(\mathcal{CCC}(G))|}\\
		&\,=\frac{1}{g(p,x)}\left(p (p^2-1)\left(p^4x^3(p^3-2p-4)+p^4x^2(p^2x-2x-5p)+x^2(2px-3p-3x-1)\right.\right. \\
		&\left.\left.+p^6x^2(x-3)+3p^5x+p^2x(5x^2-4x-6)+p^4(4x^2+3x-2)+p^3(x^3+12x^2-1)\right)\right) \\
		&\, :=\frac{f(p,x)}{g(p,x)},
	\end{align*}
	where $g(p,x)=p^4(p+2)(p^3x+p^2x-2p^2-2x-p)=p^4(p+2)(x(p^3-2)-p+p^2(x-2)) > 0$ since $p^3-2>p$ and $x>2$. Now, $f(2,x)=750x^3-1290x^2+720x-240 > 0$ for all $x \geq p^2 \geq 4$. Also, for all $p > 2$, we have $p^3-2p-4> 0, p^2x-2x-5p>0, 2px-3p-3x-1>0$ and $5x^2-4x-6>0$. Thus, $\frac{f(p,x)}{g(p,x)}>0$.
	
	If $\mathcal{CCC}(G)=(p^2+p+1)K_n$ then, by Theorem \ref{thm1}, we have
	\begin{align*}
		M_1(\mathcal{CCC}(G))&= (p^2+p+1)\frac{x(p-1)}{p^2}\left(\frac{x(p-1)}{p^2}-1\right)^2 \\
		&=\frac{x(p^3-1)}{p^6}\left(p^2x^2-2px^2+x^2-2p^3x+2p^2x+p^4\right) \text{ \quad and }
	\end{align*}
	\begin{align*}
		&M_2(\mathcal{CCC}(G))= (p^2+p+1)\frac{x(p-1)}{p^2}\frac{\left(\frac{x(p-1)}{p^2}-1\right)^3}{2} \\
		&=\frac{x(p^3-1)}{2p^8}\left(p^3x^3-3p^2x^3+3px^3-x^3-3p^4x^2+6p^3x^2-3p^2x^2+3p^5x -3p^4x-p^6\right).
	\end{align*}
	Also, $|v(\mathcal{CCC}(G))|=(p^2+p+1)\frac{x(p-1)}{p^2}=\frac{x(p^3-1)}{p^2}$ and $|e(\mathcal{CCC}(G))|=(p^2+p+1)\binom{\frac{x(p-1)}{p^2}}{2}$ $=\frac{x(p^3-1)}{2p^4}(px-x-p^2)$. 
	It can be seen that
	\[
	\frac{M_1(\mathcal{CCC}(G))}{|v(\mathcal{CCC}(G))|}=\frac{p^2x^2-2px^2+x^2-2p^3x+2p^2x+p^4}{p^4} = \frac{M_1(\mathcal{CCC}(G))}{|v(\mathcal{CCC}(G))|}.
	\]
%
	In both the cases, we have $\frac{M_{2}(\mathcal{CCC}(G))}{|e(\mathcal{CCC}(G))|} \geq \frac{M_{1}(\mathcal{CCC}(G))}{|v(\mathcal{CCC}(G))|}$ and hence the result follows.
\end{proof}
\begin{cor}
	Let $G$ be a non-abelian group of order $p^n$ such that $\frac{G}{Z(G)}$ is abelian and $|Z(G)|=p^{n-3}$, where $p$ is a prime and $n \geq 4$. 
	Then $\mathcal{CCC}(G)$ satisfies Hansen-Vuki{\v{c}}evi{\'c} conjecture. 
\end{cor}
\begin{proof}
	Clearly $|\frac{G}{Z(G)}|=p^3$ 
	and so the result follows from Theorem \ref{G/Z abelian}.
\end{proof}

In the remaining part of this section we consider non-abelian groups $G$  such that $\frac{G}{Z(G)}$ is a group of order $p^3$ for any prime  $p$. The CCC-graphs of these groups have been obtained in \cite{MAS-ARA-2020}.
\begin{theorem}\label{G/Z non-abelian}
	Let $G$ be a non-abelian group   such that $\frac{G}{Z(G)}$ is non-abelian and $|\frac{G}{Z(G)}|=p^3$, where  $p$ is a prime.
	Then $\mathcal{CCC}(G)$ satisfies Hansen-Vuki{\v{c}}evi{\'c} conjecture.  
\end{theorem}
\begin{proof}
	We have $\mathcal{CCC}(G)=K_m \sqcup kpK_{n_1} \sqcup (p-k)K_{n_2}, (kp+1)K_{n_1} \sqcup (p+1-k)K_{n_2}, K_m \sqcup pK_{n_2}, (p^2+p+1)K_{n_1}$ or $K_{n_1} \sqcup (p+1)K_{n_2}$, where $x=|Z(G)|, m=\frac{x(p^2-1)}{p}, n_1=\frac{x(p-1)}{p^2}, n_2=\frac{x(p-1)}{p}$ and $1 \leq k \leq p$. If $\mathcal{CCC}(G)=K_m \sqcup kpK_{n_1} \sqcup (p-k)K_{n_2}$ then, by Theorem \ref{thm1}, we have
		\begin{align*}
		M_1(\mathcal{CCC}(G))&=\frac{x(p^2-1)}{p}\left(\frac{x(p^2-1)}{p}-1\right)^2+kp\frac{x(p-1)}{p^2}\left(\frac{x(p-1)}{p^2}-1\right)^2 \\
		& \quad \quad \quad \quad \quad \quad \quad \quad \quad\quad \quad \quad +(p-k)\frac{x(p-1)}{p}\left(\frac{x(p-1)}{p}-1\right)^2 \\
		&=\frac{x(p-1)}{p^6}\left(p^8x^2+p^7x^2-2p^7x-p^6x^2-4p^6x-4p^5x^2+2p^4x^2+p^3x^2 \right. \\
		& \quad \quad \quad \quad \quad \quad \left. +2p^4x+p^5+4p^5x+2p^6-2kp^2x^2+kpx^2-4kp^4x \right. \\
		& \quad \quad \quad \quad \quad \quad \quad \quad \quad \quad \quad \left. +2kp^3x-kp^5x^2+2kp^4x^2+2kp^5x\right) \text{ \quad and}
	\end{align*}
	\begin{align*}
		M_2(\mathcal{CCC}(G))&=\frac{x(p^2-1)}{p}\frac{\left(\frac{x(p^2-1)}{p}-1\right)^3}{2}+kp\frac{x(p-1)}{p^2}\frac{\left(\frac{x(p-1)}{p^2}-1\right)^3}{2} \\
		& \quad \quad \quad \quad \quad \quad \quad \quad \quad\quad \quad \quad +(p-k)\frac{x(p-1)}{p}\frac{\left(\frac{x(p-1)}{p}-1\right)^3}{2} \\
		&=\frac{x(p-1)}{2p^8}\left(p^{11}x^3-3p^9x^3-2p^5x^3-3p^{10}x^2+3p^8x^2-6p^6x^2+3p^9x \right. \\
		& \quad \quad \quad \quad \quad \quad \left. -6p^7x-2p^8+p^{10}x^3-2p^8x^3+6p^6x^3-p^4x^3-3p^9x^2 \right. \\
		& \quad \quad \quad \quad \quad \quad \left. +12p^7x^2-3p^5x^2+6p^8x-3p^6x-p^7+2kp^4x^3-3kp^3x^3\right. \\
		& \quad \quad \quad \quad \quad \quad \left. +3kp^2x^3-kpx^3+6kp^4x^2-3kp^3x^2+6kp^6x-3kp^5x \right. \\
		& \quad \quad \quad \quad \quad \quad \left. -kp^7x^3+3kp^6x^3-3kp^5x^3+3kp^7x^2-6kp^6x^2-3kp^7x \right).
	\end{align*}
	Also, $|v(\mathcal{CCC}(G))|=\frac{x(p^2-1)}{p}+kp\frac{x(p-1)}{p^2}+(p-k)\frac{x(p-1)}{p}=\frac{x(p-1)}{p^2}(2p^2+p)$ and $|e(\mathcal{CCC}(G))|=\binom{\frac{x(p^2-1)}{p}}{2}+kp\binom{\frac{x(p-1)}{p^2}}{2}+(p-k)\binom{\frac{x(p-1)}{p^2}}{2}=\frac{x(p-1)}{2p^4}(p^5x+2p^4x-2p^3x-2p^4-p^2x-p^3+2p^2kx-kpx-kp^3x)$. Therefore,
	\begin{align*}
		& \frac{M_{2}(\mathcal{CCC}(G))}{|e(\mathcal{CCC}(G))|} - \frac{M_{1}(\mathcal{CCC}(G))}{|v(\mathcal{CCC}(G))|} 
		=\frac{f_1(p,k,x)}{g_1(p,k,x)},
	\end{align*}
	where
	\begin{align*}
		f_1(p,k,x)&=p^{13}x^3-4p^{11}x^3-2p^{12}x^2+2p^{10}x^2+2p^{10}x^3-2p^8x^3+2p^{11}x^2-2p^9x^2 \\
		&   -kp^6x^3-10kp^5x^3+4kp^4x^3+2kp^3x^3+8kp^6x^2+4kp^5x^2-4kp^9x^3-2kp^8x^3 \\
		&   +11kp^7x^3+4kp^9x^2+14kp^8x^2+3p^9x^3-kp^2x^3-2kp^4x^2-24kp^7x^2 \\
		&   -4kp^{10}x^2+5k^2p^4x^3-4k^2p^3x^3+12k^2p^6x^2-8k^2p^5x^2+4k^2p^7x^3 \\
		&   -5k^2p^6x^3-8k^2p^7x^2+k^2p^2x^3 +2k^2p^4x^2+kp^{11}x^3-k^2p^8x^3+2k^2p^8x^2
	\end{align*}
and $g_1(p,k,x) = p^4(2p^2+p)(p^5x+2p^4x-2p^3x-2p^4-p^2x-p^3+2p^2kx-kpx-kp^3x)$.
	On factorization, we have 
	\begin{align*}
		f_1(p,k,x)&=p^2x^2(p-1)^2\left(p^6x(p^3-p-2)+p^7(px-2p-2) +p^4x(p^4-k^2) \right. \\ 
		& \quad \quad \quad \quad \quad \quad \quad \left. +2k^2p^2(p^2+1)+kx(k-1)+kp^4(x(p^3-p-6)-4p+10) \right. \\
		& \quad \quad  \quad \quad \quad \quad \quad \quad \quad \left. +kp^2(5x-2)+2k^2p(p^2x-x-2p^2)+2kp^6(x-2)\right).
	\end{align*}
	For all prime $p \geq 2, 1 \leq k \leq p$ and $x \geq p^2$, we have $p^3-p-2 > 0, px-2p-2>0, p^4-k^2>0, p^2x-x-2p^2>0$ and $x(p^3-p-6) > 4p-10$. As such, $f_1(p,k,x) > 0$. Also, $g_1(p,k,x)=\frac{2p^{10}}{x^2(p-1)^2} \times |e(\mathcal{CCC}(G))| \times |v(\mathcal{CCC}(G))|>0$. Therefore, $\frac{f_1(p,k,x)}{g_1(p,k,x)}>0$.
	
	If $\mathcal{CCC}(G)=(kp+1)K_{n_1} \sqcup (p+1-k)K_{n_2}$ then, by Theorem \ref{thm1}, we have
	\begin{align*}
		M_1(\mathcal{CCC}(G)&=(kp+1)\frac{x(p-1)}{p^2}\left(\frac{x(p-1)}{p^2}-1\right)^2 \\
		& \quad \quad \quad \quad \quad \quad \quad \quad \quad \quad \quad \quad +(p+1-k)\frac{x(p-1)}{p}\left(\frac{x(p-1)}{p}-1\right)^2 \\
		&= \frac{x(p-1)}{p^6}\left(p^6x^2+p^6+p^5-2kp^2x^2+kpx^2-4kp^4x+2kp^3x+p^2x^2 \right. \\
		& \quad \quad \quad \quad \left. -2px^2+x^2-2p^3x+2p^2x+p^4-p^5x^2-p^4x^2-2p^6x+p^3x^2 \right. \\
		& \quad \quad \quad \quad \quad \quad \quad \quad \quad \quad \quad \left. +2p^4x-kp^5x^2+2kp^4x^2+2kp^5x \right) \text{ \quad and }
	\end{align*}
	\begin{align*}
		M_2(\mathcal{CCC}(G))&=(kp+1)\frac{x(p-1)}{p^2}\frac{\left(\frac{x(p-1)}{p^2}-1\right)^3}{2}+(p+1-k)\frac{x(p-1)}{p}\frac{\left(\frac{x(p-1)}{p}-1\right)^3}{2} \\
		&=\frac{x(p-1)}{2p^8}\left(2kp^4x^3-3kp^3x^3+3kp^2x^3-kpx^3+6kp^4x^2-3kp^3x^2+p^3x^3 \right. \\
		& \quad \quad \quad \quad \quad \left. +6kp^6x -3kp^5x-3p^2x^3+3px^3-x^3-3p^4x^2+6p^3x^2-p^6 \right. \\
		& \quad \quad \quad \quad \quad \quad \left. -3p^2x^2+3p^5x-3p^4x-2p^7x^3+2p^5x^3-3p^8x^2+3p^7x^2 \right. \\
		& \quad \quad \quad \quad \quad \quad \left. +p^8x^3+3p^8x-3p^5x^2-p^8-p^4x^3+3p^6x^2-3p^6x-p^7 \right. \\
		& \quad \quad \quad \quad \quad \quad \left. -kp^7x^3+3kp^6x^3-3kp^5x^3+3kp^7x^2-6kp^6x^2-3kp^7x \right).
	\end{align*}
	Also, $|v(\mathcal{CCC}(G))|=(kp+1)\frac{x(p-1)}{p^2}+(p+1-k)\frac{x(p-1)}{p}=\frac{x(p-1)}{p^2}(p^2+p+1)$ and $|e(\mathcal{CCC}(G))|=(kp+1)\binom{\frac{x(p-1)}{p^2}}{2}+(p+1-k)\binom{\frac{x(p-1)}{p}}{2}=\frac{x(p-1)}{2p^4}(2kp^2x-kpx+px-x-p^2+p^4x+p^4-p^2x-p^3-kp^3x)$. Therefore, 
	\begin{align*}
		&\frac{M_{2}(\mathcal{CCC}(G))}{|e(\mathcal{CCC}(G))|} - \frac{M_{1}(\mathcal{CCC}(G))}{|v(\mathcal{CCC}(G))|} 
		 =\frac{f_2(p,k,x)}{g_2(p,k,x)},
	\end{align*}
	where
	\begin{align*}
		f_2(p,k,x)&=kp^9x^3-k^2p^8x^3-3kp^8x^3+p^8x^3+4k^2p^7x^3-3p^7x^3-5k^2p^6x^3+9kp^6x^3 \\
		& \quad +p^6x^3-10kp^5x^3+5p^5x^3+5k^2p^4x^3-kp^4x^3-5p^4x^3-4k^2p^3x^3+8kp^3x^3 \\
		& \quad -p^3x^3+k^2p^2x^3-5kp^2x^3+3p^2x^3+kpx^3-px^3-2kp^9x^2+2k^2p^8x^2 \\
		& \quad +6kp^8x^2-2p^8x^2-8k^2p^7x^2-2kp^7x^2+6p^7x^2+12k^2p^6x^2-12kp^6x^2\\
		& \quad -4p^6x^2-8k^2p^5x^2+18kp^5x^2 -4p^5x^2+2k^2p^4x^2-10kp^4x^2+6p^4x^2 \\
		& \quad \quad \quad \quad \quad \quad \quad \quad \quad \quad \quad \quad \quad \quad \quad \quad \quad \quad \quad \quad \quad \quad \quad +2kp^3x^2-2p^3x^2
	\end{align*}
and $g_2(p,k,x) = p^4(p^2+p+1)(2kp^2x-kpx+px-x-p^2+p^4x+p^4-p^2x-p^3-kp^3x)$.
	On factorization, we have $f_2(p,k,x)=px^2(p-1)^4(p+1-k)(kp+1)(p^2x-x-2p^2)$. For all $p \geq 2, 1 \leq k \leq p$ and $x \geq p^2$, we have $p+1-k > 0$ and $x(p^2-1) > 2p^2$. As such, $f_2(p,k,x) > 0$. Also, $g_2(p,k,x)=\frac{2p^{10}}{x^2(p-1)^2} \times |e(\mathcal{CCC}(G))| \times |v(\mathcal{CCC}(G))|>0$. Therefore, $\frac{f_2(p,k,x)}{g_2(p,k,x)}>0$.
	
	If $\mathcal{CCC}(G)=K_m \sqcup pK_{n_2}$ then, by Theorem \ref{thm1}, we have
	\begin{align*}
		M_1(\mathcal{CCC}(G))&=\frac{x(p^2-1)}{p}\left(\frac{x(p^2-1)}{p}-1\right)^2+p \frac{x(p-1)}{p}\left(\frac{x(p-1)}{p}-1\right)^2 \\
		&=\frac{x(p-1)}{p^3}\left(p^5x^2-p^3x^2+2px^2-2p^4x+2p^3+p^4x^2-4p^2x^2+x^2 \right. \\
		& \quad \quad \quad \quad \quad \quad \quad \quad \quad \quad \quad \quad \quad \left. +2px-4p^3x+4p^2x+p^2 \right)  
	\end{align*} 
and
	\begin{align*}
		M_2&(\mathcal{CCC}(G))=\frac{x(p^2-1)}{p}\frac{\left(\frac{x(p^2-1)}{p}-1\right)^3}{2}+p \frac{x(p-1)}{p}\frac{\left(\frac{x(p-1)}{p}-1\right)^3}{2} \\
		&=\frac{x(p-1)}{2p^4}\left(p^7x^3-3p^5x^3-2px^3-3p^6x^2+3p^4x^2-6p^2x^2+3p^5x-6p^3x \right. \\
		&   \left. -2p^4+p^6x^3-2p^4x^3+6p^2x^3-x^3-3p^5x^2+12p^3x^2-3px^2 +6p^4x-3p^2x-p^3\right. \\
		& \quad \quad \quad \quad \quad \quad  \quad \quad \quad \quad \quad \quad  \quad \quad \quad \quad \quad \quad  \quad \quad \quad \left. +6p^4x-3p^2x-p^3\right).
	\end{align*}
	Also, $|v(\mathcal{CCC}(G))|=\frac{x(p^2-1)}{p}+p\frac{x(p-1)}{p}=\frac{x(p-1)}{p}(2p+1)$ and $|e(\mathcal{CCC}(G))|=\binom{\frac{x(p^2-1)}{p}}{2}+p\binom{\frac{x(p-1)}{p}}{2}=\frac{x(p-1)}{2p^2}(p^3x-2px+2p^2x-2p^2-x-p)$. Therefore, for all $p \geq 2$ and $x \geq p^2$, we have
	\begin{align*}
		&\frac{M_{2}(\mathcal{CCC}(G))}{|e(\mathcal{CCC}(G))|} - \frac{M_{1}(\mathcal{CCC}(G))}{|v(\mathcal{CCC}(G))|} \\
		&\qquad\qquad=\frac{p^6x^3(p^2-4)+2p^5x^2(x-p^2)+p^3x^3(3p-2)+2p^6x^2+2p^4x^2(p-1)}{p^2(2p+1)\{x(p^3-2p-1)+p(2px-2p-1)\}}  >0.
	\end{align*}
	If $\mathcal{CCC}(G)=(p^2+p+1)K_{n_1}$ then the result follows from Theorem \ref{G/Z abelian}.
	If $\mathcal{CCC}(G)=K_{n_1} \sqcup (p+1)K_{n_2}$ then, by Theorem \ref{thm1}, we have
	\begin{align*}
		M_1(\mathcal{CCC}(G))&=\frac{x(p-1)}{p^2}\left(\frac{x(p-1)}{p^2}-1\right)^2+(p+1)\frac{x(p-1)}{p}\left(\frac{x(p-1)}{p}-1\right)^2 \\
		&=\frac{x(p-1)}{p^6}\left(p^6x^2-p^5x^2-p^4x^2+p^3x^2-2p^6x+p^6+2p^4x+p^5+p^4 \right. \\
		& \quad \quad \quad \quad \quad \quad \quad \quad \quad \quad \quad \left. +p^2x^2-2px^2+x^2-2p^3x+2p^2x \right)
	\end{align*}
and

	\begin{align*}
		M_2(\mathcal{CCC}(G))&=\frac{x(p-1)}{p^2}\frac{\left(\frac{x(p-1)}{p^2}-1\right)^3}{2}+(p+1)\frac{x(p-1)}{p}\frac{\left(\frac{x(p-1)}{p}-1\right)^3}{2} \\
		&= \frac{x(p-1)}{2p^8}\left(p^8x^3-2p^7x^3+2p^5x^3-p^4x^3-3p^8x^2+3p^7x^2+3p^6x^2 \right. \\
		& \quad \quad \quad \quad \quad \left. -3p^5x^2+3p^8x-3p^6x+3p^5x-3p^4x-p^8-p^7-p^6+p^3x^3 \right. \\
		& \quad \quad \quad \quad \quad \quad \quad \quad \quad \quad \left. -3p^2x^3+3px^3-x^3-3p^4x^2+6p^3x^2-3p^2x^2 \right).
	\end{align*}
	Also, $|v(\mathcal{CCC}(G))|=\frac{x(p-1)}{p^2}+(p+1)\frac{x(p-1)}{p}=\frac{x(p-1)}{p^2}(p^2+p+1)$ and $|e(\mathcal{CCC}(G))|=\binom{\frac{x(p-1)}{p^2}}{2}+(p+1)\binom{\frac{x(p-1)}{p}}{2}=\frac{x(p-1)}{2p^4}(p^4x+px-p^2x-p^4-p^3-p^2-x)$. Therefore, $\frac{M_{2}(\mathcal{CCC}(G))}{|e(\mathcal{CCC}(G))|} - \frac{M_{1}(\mathcal{CCC}(G))}{|v(\mathcal{CCC}(G))|}=\frac{f_3(p,x)}{p^4(p^2+p+1)(p^4x+px-p^2x-p^4-p^3-p^2-x)}:=\frac{f_3(p,x)}{g_3(p,x)}$, where

	\begin{align*}
		f_3(p,x)&=p^8x^3-3p^7x^3+p^6x^3+5p^5x^3-5p^4x^3-p^3x^3+3p^2x^3-px^3-2p^8x^2 \\
		& \quad \quad \quad \quad \quad \quad \quad \quad \quad \quad +6p^7x^2-4p^6x^2-4p^5x^2+6p^4x^2-2p^3x^2.
	\end{align*}
	On factorization, we have $f_3(p,x)=px^2(p+1)(p-1)^4(p^2x-x-2p^2)$. For all prime $p \geq 2$ and $x \geq p^2$, we have $x(p^2-1)>2p^2$. As such, $f_3(p,x)> 0$ and $g_3(p,x)=\frac{2p^{10}}{x^2(p-1)^2} \times |e(\mathcal{CCC}(G))| \times |v(\mathcal{CCC}(G))|>0$. Therefore, $\frac{f_3(p,x)}{g_3(p,x)}>0$.
	
	In all the cases, we have $\frac{M_{2}(\mathcal{CCC}(G))}{|e(\mathcal{CCC}(G))|} > \frac{M_{1}(\mathcal{CCC}(G))}{|v(\mathcal{CCC}(G))|}$. Hence, the result follows.
\end{proof}
\begin{cor}\label{cor G/Z non-abelian}
	Let $G$ be a non-abelian group of order $p^n$ such that $\frac{G}{Z(G)}$ is non-abelian  and $|Z(G)|=p^{n-3}$, where $p$ is a prime and $n \geq 4$. 
	Then $\mathcal{CCC}(G)$ satisfies Hansen-Vuki{\v{c}}evi{\'c} conjecture. 
\end{cor}
\begin{proof}
	Clearly $|\frac{G}{Z(G)}|=p^3$ and so the result follows from Theorem \ref{G/Z non-abelian}.
\end{proof}
\begin{cor}
	If $G$ is a non-abelian group of order $p^4$ then 
	$\mathcal{CCC}(G)$ satisfies Hansen-Vuki{\v{c}}evi{\'c} conjecture.  
\end{cor}
\begin{proof}
	We have $|G|=p^4$, $p$ is a prime, then $|Z(G)|$ must be $p$ or $p^2$. If $|Z(G)|=p^2$ then the result follows from Corollary \ref{order G/Z=p^2}. Otherwise, if $|Z(G)|=p$ then the result follows from Corollary \ref{cor G/Z non-abelian}.
\end{proof}

\section{Concluding remarks}

In this paper we have shown that $\mathcal{CCC}(G)$ satisfies Hansen-Vuki{\v{c}}evi{\'c} conjecture if 
	\begin{itemize}
		\item $G$ is isomorphic to $D_{2m}, Q_{4m}, SD_{8m}, V_{8m}, U(n,m)$ and $G(p,m,n)$.
		\item $\frac{G}{Z(G)}$ is isomorphic to $D_{2m}, \mathbb{Z}_p \times \mathbb{Z}_p$ and Frobenious group of order $pq$ and $p^2q$ for any two primes $p$ and $q$.
		\item $|\frac{G}{Z(G)}|=p^3$ for any prime $p.$ 
	\end{itemize}
Moreover, it was also seen that if $G$ is any non-abelian group of order $p^3$ or $p^4$, where $p$ is any prime, then $\mathcal{CCC}(G)$  satisfies Hansen-Vuki{\v{c}}evi{\'c} conjecture. Using  Theorem \ref{cen-iso-D}, Theorem \ref{G/Z=Z_p x Z_p} and the characterizations of $4$-centralizer (``$G$ is a $4$-centralizer finite group if and only if $\frac{G}{Z(G)} \cong \mathbb{Z}_2 \times \mathbb{Z}_2$" \cite[Theorem 2]{BS-1994}),  $5$-centralizer (``$G$ is a $5$-centralizer finite group if and only if $\frac{G}{Z(G)} \cong \mathbb{Z}_3 \times \mathbb{Z}_3$ or $D_{2\times 3}$" \cite[Theorem 4]{BS-1994})  finite groups we have  the following result.
\begin{theorem}
If $G$ is a $4$-centralizer or $5$-centralizer finite group then $\mathcal{CCC}(G)$  satisfies Hansen-Vuki{\v{c}}evi{\'c} conjecture.
\end{theorem}
Using Theorem \ref{G/Z=Z_p x Z_p} and the characterization of $(p+2)$-centralizer finite $p$-group (``If $G$ is a finite $p$-group then it is $(p+2)$-centralizer if and only if $\frac{G}{Z(G)} \cong \mathbb{Z}_p \times \mathbb{Z}_p$" \cite[Lemma 2.7]{Ashrafi-2000})   we have  the following result.
\begin{theorem}
	If $G$ is a $(p+2)$-centralizer finite $p$-group then $\mathcal{CCC}(G)$  satisfies Hansen-Vuki{\v{c}}evi{\'c} conjecture.
\end{theorem}
If $\{g_1, g_2, \dots, g_m\}$ is a set of pairwise non-commuting elements of $G$ having maximal size then $G$ is a $4$-centralizer group if $m = 3$ and a $5$-centralizer group if $m = 4$ (see \cite[Lemma 2.4]{AJH-2007}). Therefore, if $m = 3$ and $4$ then $\mathcal{CCC}(G)$ graphs of such groups satisfy Hansen-Vuki{\v{c}}evi{\'c} conjecture.

Let $p$ be the smallest prime divisor of $|G|$. Then the fact that ``commuting probability of $G$, $\Pr(G) = \frac{p^2 + p - 1}{p^3}$ if and only if $\frac{G}{Z(G)} \cong \mathbb{Z}_p \times \mathbb{Z}_p$" \cite[Theorem 3]{MacHale-1974} along with Theorem \ref{G/Z=Z_p x Z_p}  give the following result.
\begin{theorem}
If $p$ is the smallest prime divisor of $|G|$ and  $\Pr(G) = \frac{p^2 + p - 1}{p^3}$  then $\mathcal{CCC}(G)$  satisfies Hansen-Vuki{\v{c}}evi{\'c} conjecture.
\end{theorem}
In addition, if $\Pr(G) = \frac{5}{14}, \frac{2}{5}, \frac{11}{27},  \frac{1}{2}$ or  $\frac{5}{8}$ then  $\frac{G}{Z(G)} \cong   D_{2\times 7}, D_{2\times 5}, D_{2\times 4}, D_{2\times 3}$ or  $\mathbb{Z}_2 \times \mathbb{Z}_2$ (see \cite[pp. 246]{Rusin-1979} and \cite[pp. 451]{Nath-2013}). Therefore, the following result follows in view of Theorem \ref{cen-iso-D} and Theorem \ref{G/Z=Z_p x Z_p}.
\begin{theorem}
If $\Pr(G) = \frac{5}{14}, \frac{2}{5}, \frac{11}{27},  \frac{1}{2}$ or  $\frac{5}{8}$  then $\mathcal{CCC}(G)$  satisfies Hansen-Vuki{\v{c}}evi{\'c} conjecture.
\end{theorem}

It is worth noting that in this work we have not found any group $G$ such that $\mathcal{CCC}(G)$ does not satisfy Hansen-Vuki{\v{c}}evi{\'c} conjecture. We conclude this paper with the following question.

\begin{question}
Is there any finite group $G$ such that $\mathcal{CCC}(G)$ does not satisfy Hansen-Vuki{\v{c}}evi{\'c} conjecture?
\end{question}

\vspace{.5cm}

\noindent {\bf Acknowledgement.} The first author gratefully acknowledges the Council of Scientific and Industrial Research for awarding the fellowship (File No. 09/0796(16521)/2023-EMR-I).

\noindent {\bf Conflict of interest.} The authors declare that they have no conflict of interest.

\noindent {\bf Funding information.} No funding was received by the authors.


\end{document}